\documentclass{amsart}
\usepackage{amssymb,amscd,graphics,axodraw}

\numberwithin{equation}{section}

\newtheorem{prop}{Proposition}
\newtheorem{theorem}[prop]{Theorem}

\newtheorem{lemma}[prop]{Lemma}

\theoremstyle{definition}

\newtheorem{remark}[prop]{Remark}
\newtheorem{assumption}[prop]{Assumption}

\numberwithin{prop}{section}

\newcommand{\comp}{\mathrm{comp}}
\newcommand{\f}{\upsilon}
\newcommand{\la}{\lambda}
\newcommand{\La}{\Lambda}
\newcommand{\lb}{\bar{\ell}}

\newcommand{\lt}{\tilde{\ell}}
\newcommand{\mut}{\tilde{\mu}}
\newcommand{\nt}{\nu^\bullet}
\newcommand{\ntt}{{\tilde{\nu}}^\bullet}
\newcommand{\qbin}[2]{\genfrac{[}{]}{0pt}{}{#1}{#2}}
\newcommand{\rk}{\mathrm{rk}}
\newcommand{\wt}{\mathrm{wt}}

\newcommand{\I}{I^*}
\newcommand{\Jm}{J_{\mathrm{max}}}
\newcommand{\Jt}{J^{\bullet}}
\newcommand{\Jtt}{{\tilde{J}}^{\bullet}}
\newcommand{\N}{\Z_{\ge0}}
\newcommand{\Path}{\mathcal{P}}
\newcommand{\RC}{\mathrm{RC}}
\newcommand{\Z}{\mathbb{Z}}

\newcommand{\geh}{\mathfrak{g}}
\newcommand{\gehb}{\overline{\geh}}
\newcommand{\gf}{\geh_{\overline{0}}}
\newcommand{\Lab}{\overline{\Lambda}}
\newcommand{\Phit}{\tilde{\Phi}}
\newcommand{\Pfin}{\overline{P}}

\newcommand{\vphi}{\varphi}
\newcommand{\Xb}{\overline{X}}
\newcommand{\Mb}{\overline{M}}

\newcommand{\Db}{\overline{D}}
\newcommand{\Hb}{\overline{H}}

\newcommand{\alt}{\tilde{\alpha}}
\newcommand{\Lt}{\tilde{\La}}
\newcommand{\Pt}{\tilde{P}}

\newcommand{\cd}{\dotsm}

\newcommand{\et}[1]{\widetilde{e}_{#1}}
\newcommand{\ft}[1]{\widetilde{f}_{#1}}

\newcommand{\HHH}{\mathcal{H}}

\begin{document}

\title{A crystal to rigged configuration bijection for nonexceptional affine
algebras}

\author[M.~Okado]{Masato Okado}
\address{Department of Informatics and Mathematical Science,
Graduate School of Engineering Science, Osaka University,
Toyonaka, Osaka 560-8531, Japan}
\email{okado@sigmath.es.osaka-u.ac.jp}

\author[A.~Schilling]{Anne Schilling}
\address{Department of Mathematics, University of California, One Shields
Avenue, Davis, CA 95616-8633, U.S.A.}
\email{anne@math.ucdavis.edu}

\author[M.~Shimozono]{Mark Shimozono}
\address{Department of Mathematics, 460 McBryde Hall, Virginia Tech,
Blacksburg, VA 24061-0123, U.S.A}
\email{mshimo@math.vt.edu}


\begin{abstract}
Kerov, Kirillov, and Reshetikhin defined a bijection between
highest weight vectors in the crystal graph of a tensor power of
the vector representation, and combinatorial objects called rigged
configurations, for type $A^{(1)}_n$. We define an analogous
bijection for all nonexceptional affine types, thereby proving (in
this special case) the fermionic formulas conjectured by Hatayama,
Kuniba, Takagi, Tsuboi, Yamada, and the first author.
\end{abstract}

\maketitle


\section{Introduction}
The fermionic formula, denoted by $M$, is a certain polynomial
expressed as a sum of products of $q$-binomial coefficients. It
originates in the Bethe Ansatz analysis of solvable lattice models
in two dimensional statistical mechanics. The prototypical example
is given by the Kostka polynomial $K_{\la\mu}(q)\in\Z_{\ge0}[q]$,
which is indexed by a pair of partitions $\la,\mu$. According to
Lascoux and Sch\"utzenberger \cite{LS},
\begin{equation*}
K_{\la\mu}(q)=\sum_{\la\in\mathcal{T}(\la,\mu)}q^{c(T)}.
\end{equation*}
Here $\mathcal{T}(\la,\mu)$ is the set of semistandard tableaux of
shape $\la$ and weight $\mu$, and $c(T)$ is the charge of the
tableau $T$.

We consider the case that $\mu$ is a single column $(1^L)$.
Kirillov and Reshetikhin \cite{KR} gave a fermionic formula for
the Kostka polynomial:
\begin{equation} \label{eq:KM}
K_{\la,(1^L)}(q) = q^{\binom{L}{2}} M(\la,(1^L);q^{-1})
\end{equation}
where
\begin{equation*}
\begin{split}
M(\la,(1^L);q) &=\sum_{\{m\}}q^{cc(\{m\})}\prod_{\substack{1\le
a\le n\\ i\ge1}} \begin{bmatrix} p^{(a)}_i+m^{(a)}_i \\
m^{(a)}_i \end{bmatrix}, \\
cc(\{m\})&=\frac{1}{2}\sum_{1\le a,b\le n}C_{ab}
\sum_{i,j\ge1}\min(i,j)m^{(a)}_im^{(b)}_j, \\
p^{(a)}_i&=L\delta_{a1}-\sum_{1\le b\le
n}C_{ab}\sum_{j\ge1}\min(i,j)m^{(b)}_j,
\end{split}
\end{equation*}
$\begin{bmatrix} p+m \\ m \end{bmatrix}=(q)_{p+m}/(q)_p(q)_m$ is the
$q$-binomial coefficient, $(q)_m=(1-q)(1-q^2)\cdots(1-q^m)$, the sum
$\sum_{\{m\}}$ is taken
over $\{m^{(a)}_i\in\Z_{\ge0}\mid 1\le a\le n,i\ge1\}$, satisfying
$p^{(a)}_i\ge0$ for $1\le a\le n,i\ge1$ and
$\sum_{i\ge1}im^{(a)}_i=\la_{a+1}+\la_{a+2}+\cdots+\la_{n+1}$ for
$1\le a\le n$. Here $n$ is an integer not less than the length of
$\la$ minus 1, and $(C_{ab})_{1\le a,b\le n}$ is the Cartan matrix
of $sl_{n+1}$.

To prove that the Kostka polynomial is given by the fermionic
formula, Kerov, Kirillov and Reshetikhin (KKR) defined a bijection
between $\mathcal{T}(\la,(1^L))$ and combinatorial objects called
rigged configurations \cite{KKR}. Expanding the $q$-binomial
coefficients in $M(\la,(1^L);q)$, to each term $q^c$ one can
associate a rigged configuration having a statistic $c$.
Under the bijection, the charge of a tableau agrees with the
statistic on the rigged configuration. This bijection was extended
to the larger class of Littlewood-Richardson tableaux and
corresponding rigged configurations \cite{KSS}.

The Kostka polynomial is related to the affine Lie algebra of type
$A^{(1)}_n$, since the corresponding fermionic formula is derived
from the integrable model associated to the quantum affine algebra
$U_q(A^{(1)}_n)$. The Kostka polynomial $K_{\la\mu}(q)$ gives the
graded multiplicity of the $\la$-th irreducible $U_q(A_n)$-module
in the restriction of the tensor product of certain
finite-dimensional $U'_q(A^{(1)}_n)$-modules that have crystal
bases. The situation generalizes to the context of any affine Lie
algebra. One can define the analogous tensor product modules and
graded multiplicities, and a corresponding fermionic formula $M$
\cite{HKOTT,HKOTY}. The new combinatorial objects which replace
tableaux are called paths. A path is a highest weight element of
the aforementioned tensor product crystal base. Paths have a
natural statistic called energy. In the case of the Kostka
polynomial, paths biject with rigged configurations: one may send
the path (which may be viewed as a word) to its Robinson-Schensted
recording tableau, which is then sent to a rigged configuration by
the KKR bijection. The generating function of paths by energy is
called the ``one dimensional sum" $X$. The equality $X=M$ was
conjectured in full generality in \cite{HKOTT,HKOTY}.

The purpose of the paper is to construct the analogue of the KKR
bijection and thereby prove the $X=M$ conjecture, for all
nonexceptional affine Lie algebras, in the case of the simplest
crystal bases. For $A^{(1)}_n$, this case corresponds to the
Kostka polynomial $K_{\la(1^L)}(q)$ discussed above.

\section{Quantum affine algebras and crystals}
\subsection{Affine algebras}
\label{subsec:algebra} We adopt the notation of \cite{HKOTT}. Let
$\geh$ be a Kac-Moody Lie algebra of nonexceptional affine type
$X^{(r)}_N$, that is, one of the types $A^{(1)}_n (n \ge 1)$,
$B^{(1)}_n (n \ge 3)$, $C^{(1)}_n (n \ge 2)$, $D^{(1)}_n (n \ge
4)$, $A^{(2)}_{2n}(n\ge1)$, $A^{(2)\dagger}_{2n}(n\ge1)$,
$A^{(2)}_{2n-1}(n\ge2)$, $D^{(2)}_{n+1}(n\ge2)$. Note that
$A^{(2)\dagger}_{2n}$ is the same diagram as $A^{(2)}_{2n}$ but
with the opposite labeling.

The Dynkin diagram of $\geh = X^{(r)}_N$ is depicted in Table
\ref{tab:Dynkin} (Table Aff 1-3 in \cite{Kac}). Its nodes are
labeled by the set $I=\{0,1,2\dotsc,n\}$.

{\unitlength=.95pt
\begin{table}
\caption{Dynkin diagrams for $X^{(r)}_N$. The labeling of the
nodes (by elements of $I$) is specified under or the right side of
the nodes. The numbers $t_i$ (resp. $t^\vee_i$) defined in
\eqref{eq:ttdef} are attached \textit{above} the nodes for $r=1$
(resp. $r>1$) if and only if $t_i \neq 1$ (resp. $t^\vee_i \neq
1$).} \label{tab:Dynkin}
\begin{tabular}[t]{rl}
$A_1^{(1)}$:&
\begin{picture}(26,20)(-5,-5)
\multiput( 0,0)(20,0){2}{\circle{6}}
\multiput(2.85,-1)(0,2){2}{\line(1,0){14.3}}
\put(0,-5){\makebox(0,0)[t]{$0$}}
\put(20,-5){\makebox(0,0)[t]{$1$}} \put( 6, 0){\makebox(0,0){$<$}}
\put(14, 0){\makebox(0,0){$>$}}
\end{picture}
\\
&
\\
\begin{minipage}{4em}
\begin{flushright}
$A_n^{(1)}$:\\$(n \ge 2)$
\end{flushright}
\end{minipage}&
\begin{picture}(106,40)(-5,-5)
\multiput( 0,0)(20,0){2}{\circle{6}}
\multiput(80,0)(20,0){2}{\circle{6}} \put(50,20){\circle{6}}
\multiput( 3,0)(20,0){2}{\line(1,0){14}}
\multiput(63,0)(20,0){2}{\line(1,0){14}}
\multiput(39,0)(4,0){6}{\line(1,0){2}}
\put(2.78543,1.1142){\line(5,2){44.429}}
\put(52.78543,18.8858){\line(5,-2){44.429}}
\put(0,-5){\makebox(0,0)[t]{$1$}}
\put(20,-5){\makebox(0,0)[t]{$2$}}
\put(80,-5){\makebox(0,0)[t]{$n\!\! -\!\! 1$}}
\put(100,-5){\makebox(0,0)[t]{$n$}}
\put(55,20){\makebox(0,0)[lb]{$0$}}
\end{picture}
\\
&
\\
\begin{minipage}[b]{4em}
\begin{flushright}
$B_n^{(1)}$:\\$(n \ge 3)$
\end{flushright}
\end{minipage}&
\begin{picture}(126,40)(-5,-5)
\multiput( 0,0)(20,0){3}{\circle{6}}
\multiput(100,0)(20,0){2}{\circle{6}} \put(20,20){\circle{6}}
\multiput( 3,0)(20,0){3}{\line(1,0){14}}
\multiput(83,0)(20,0){1}{\line(1,0){14}}
\put(20,3){\line(0,1){14}}
\multiput(102.85,-1)(0,2){2}{\line(1,0){14.3}} 
\multiput(59,0)(4,0){6}{\line(1,0){2}} 
\put(110,0){\makebox(0,0){$>$}} \put(0,-5){\makebox(0,0)[t]{$1$}}
\put(20,-5){\makebox(0,0)[t]{$2$}}
\put(40,-5){\makebox(0,0)[t]{$3$}}
\put(100,-5){\makebox(0,0)[t]{$n\!\! -\!\! 1$}}
\put(120,-5){\makebox(0,0)[t]{$n$}}
\put(25,20){\makebox(0,0)[l]{$0$}}

\put(120,13){\makebox(0,0)[t]{$2$}}
\end{picture}
\\
&
\\
\begin{minipage}[b]{4em}
\begin{flushright}
$C_n^{(1)}$:\\$(n \ge 2)$
\end{flushright}
\end{minipage}&
\begin{picture}(126,20)(-5,-5)
\multiput( 0,0)(20,0){3}{\circle{6}}
\multiput(100,0)(20,0){2}{\circle{6}}
\multiput(23,0)(20,0){2}{\line(1,0){14}}
\put(83,0){\line(1,0){14}}
\multiput( 2.85,-1)(0,2){2}{\line(1,0){14.3}} 
\multiput(102.85,-1)(0,2){2}{\line(1,0){14.3}} 
\multiput(59,0)(4,0){6}{\line(1,0){2}} 
\put(10,0){\makebox(0,0){$>$}} \put(110,0){\makebox(0,0){$<$}}
\put(0,-5){\makebox(0,0)[t]{$0$}}
\put(20,-5){\makebox(0,0)[t]{$1$}}
\put(40,-5){\makebox(0,0)[t]{$2$}}
\put(100,-5){\makebox(0,0)[t]{$n\!\! -\!\! 1$}}
\put(120,-5){\makebox(0,0)[t]{$n$}}

\put(20,13){\makebox(0,0)[t]{$2$}}
\put(40,13){\makebox(0,0)[t]{$2$}}
\put(100,13){\makebox(0,0)[t]{$2$}}
\end{picture}
\\
&
\\
\begin{minipage}[b]{4em}
\begin{flushright}
$D_n^{(1)}$:\\$(n \ge 4)$
\end{flushright}
\end{minipage}&
\begin{picture}(106,40)(-5,-5)
\multiput( 0,0)(20,0){2}{\circle{6}}
\multiput(80,0)(20,0){2}{\circle{6}}
\multiput(20,20)(60,0){2}{\circle{6}} \multiput(
3,0)(20,0){2}{\line(1,0){14}}
\multiput(63,0)(20,0){2}{\line(1,0){14}}
\multiput(39,0)(4,0){6}{\line(1,0){2}}
\multiput(20,3)(60,0){2}{\line(0,1){14}}
\put(0,-5){\makebox(0,0)[t]{$1$}}
\put(20,-5){\makebox(0,0)[t]{$2$}}
\put(80,-5){\makebox(0,0)[t]{$n\!\! - \!\! 2$}}
\put(103,-5){\makebox(0,0)[t]{$n\!\! -\!\! 1$}}
\put(25,20){\makebox(0,0)[l]{$0$}}
\put(85,20){\makebox(0,0)[l]{$n$}}
\end{picture}
\\
&
\\

$A^{(2)}_2$:&
\begin{picture}(26,20)(-5,-5)
\multiput( 0,0)(20,0){2}{\circle{6}}
\multiput(2.958,-0.5)(0,1){2}{\line(1,0){14.084}}
\multiput(2.598,-1.5)(0,3){2}{\line(1,0){14.804}}
\put(0,-5){\makebox(0,0)[t]{$0$}}
\put(20,-5){\makebox(0,0)[t]{$1$}} \put(10,0){\makebox(0,0){$<$}}
\put(0,13){\makebox(0,0)[t]{$2$}}
\put(20,13){\makebox(0,0)[t]{$2$}}
\end{picture}
\\
&
\\
\begin{minipage}[b]{4em}
\begin{flushright}
$A_{2n}^{(2)}$:\\$(n \ge 2)$
\end{flushright}
\end{minipage}&
\begin{picture}(126,20)(-5,-5)
\multiput( 0,0)(20,0){3}{\circle{6}}
\multiput(100,0)(20,0){2}{\circle{6}}
\multiput(23,0)(20,0){2}{\line(1,0){14}}
\put(83,0){\line(1,0){14}}
\multiput( 2.85,-1)(0,2){2}{\line(1,0){14.3}} 
\multiput(102.85,-1)(0,2){2}{\line(1,0){14.3}} 
\multiput(59,0)(4,0){6}{\line(1,0){2}} 
\put(10,0){\makebox(0,0){$<$}} \put(110,0){\makebox(0,0){$<$}}
\put(0,-5){\makebox(0,0)[t]{$0$}}
\put(20,-5){\makebox(0,0)[t]{$1$}}
\put(40,-5){\makebox(0,0)[t]{$2$}}
\put(100,-5){\makebox(0,0)[t]{$n\!\! -\!\! 1$}}
\put(120,-5){\makebox(0,0)[t]{$n$}}
\put(0,13){\makebox(0,0)[t]{$2$}}
\put(20,13){\makebox(0,0)[t]{$2$}}
\put(40,13){\makebox(0,0)[t]{$2$}}
\put(100,13){\makebox(0,0)[t]{$2$}}
\put(120,13){\makebox(0,0)[t]{$2$}}
\put(120,13){\makebox(0,0)[t]{$2$}}
\end{picture}
\\
&
\\
$A^{(2)\dagger}_2$:&
\begin{picture}(26,20)(-5,-5)
\multiput( 0,0)(20,0){2}{\circle{6}}
\multiput(2.958,-0.5)(0,1){2}{\line(1,0){14.084}}
\multiput(2.598,-1.5)(0,3){2}{\line(1,0){14.804}}
\put(0,-5){\makebox(0,0)[t]{$0$}}
\put(20,-5){\makebox(0,0)[t]{$1$}} \put(10,0){\makebox(0,0){$>$}}
\put(0,13){\makebox(0,0)[t]{$$}} \put(20,13){\makebox(0,0)[t]{$$}}
\end{picture}
\\
&
\\
\begin{minipage}[b]{4em}
\begin{flushright}
$A_{2n}^{(2)\dagger}$:\\$(n \ge 2)$
\end{flushright}
\end{minipage}&
\begin{picture}(126,20)(-5,-5)
\multiput( 0,0)(20,0){3}{\circle{6}}
\multiput(100,0)(20,0){2}{\circle{6}}
\multiput(23,0)(20,0){2}{\line(1,0){14}}
\put(83,0){\line(1,0){14}}
\multiput( 2.85,-1)(0,2){2}{\line(1,0){14.3}} 
\multiput(102.85,-1)(0,2){2}{\line(1,0){14.3}} 
\multiput(59,0)(4,0){6}{\line(1,0){2}} 
\put(10,0){\makebox(0,0){$>$}} \put(110,0){\makebox(0,0){$>$}}
\put(0,-5){\makebox(0,0)[t]{$0$}}
\put(20,-5){\makebox(0,0)[t]{$1$}}
\put(40,-5){\makebox(0,0)[t]{$2$}}
\put(100,-5){\makebox(0,0)[t]{$n\!\! -\!\! 1$}}
\put(120,-5){\makebox(0,0)[t]{$n$}}
\put(0,13){\makebox(0,0)[t]{$$}}
\put(20,13){\makebox(0,0)[t]{$$}}
\put(40,13){\makebox(0,0)[t]{$$}}
\put(100,13){\makebox(0,0)[t]{$$}}
\put(120,13){\makebox(0,0)[t]{$$}}
\put(120,13){\makebox(0,0)[t]{$$}}
\end{picture}
\\

\begin{minipage}[b]{4em}
\begin{flushright}
$A_{2n-1}^{(2)}$:\\$(n \ge 3)$
\end{flushright}
\end{minipage}&
\begin{picture}(126,40)(-5,-5)
\multiput( 0,0)(20,0){3}{\circle{6}}
\multiput(100,0)(20,0){2}{\circle{6}} \put(20,20){\circle{6}}
\multiput( 3,0)(20,0){3}{\line(1,0){14}}
\multiput(83,0)(20,0){1}{\line(1,0){14}}
\put(20,3){\line(0,1){14}}
\multiput(102.85,-1)(0,2){2}{\line(1,0){14.3}} 
\multiput(59,0)(4,0){6}{\line(1,0){2}} 
\put(110,0){\makebox(0,0){$<$}} \put(0,-5){\makebox(0,0)[t]{$1$}}
\put(20,-5){\makebox(0,0)[t]{$2$}}
\put(40,-5){\makebox(0,0)[t]{$3$}}
\put(100,-5){\makebox(0,0)[t]{$n\!\! -\!\! 1$}}
\put(120,-5){\makebox(0,0)[t]{$n$}}
\put(25,20){\makebox(0,0)[l]{$0$}}

\put(120,13){\makebox(0,0)[t]{$2$}}
\end{picture}
\\
&
\\
\begin{minipage}[b]{4em}
\begin{flushright}
$D_{n+1}^{(2)}$:\\$(n \ge 2)$
\end{flushright}
\end{minipage}&
\begin{picture}(126,20)(-5,-5)
\multiput( 0,0)(20,0){3}{\circle{6}}
\multiput(100,0)(20,0){2}{\circle{6}}
\multiput(23,0)(20,0){2}{\line(1,0){14}}
\put(83,0){\line(1,0){14}}
\multiput( 2.85,-1)(0,2){2}{\line(1,0){14.3}} 
\multiput(102.85,-1)(0,2){2}{\line(1,0){14.3}} 
\multiput(59,0)(4,0){6}{\line(1,0){2}} 
\put(10,0){\makebox(0,0){$<$}} \put(110,0){\makebox(0,0){$>$}}
\put(0,-5){\makebox(0,0)[t]{$0$}}
\put(20,-5){\makebox(0,0)[t]{$1$}}
\put(40,-5){\makebox(0,0)[t]{$2$}}
\put(100,-5){\makebox(0,0)[t]{$n\!\! -\!\! 1$}}
\put(120,-5){\makebox(0,0)[t]{$n$}}

\put(20,13){\makebox(0,0)[t]{$2$}}
\put(40,13){\makebox(0,0)[t]{$2$}}
\put(100,13){\makebox(0,0)[t]{$2$}}
\end{picture}
\\
&
\\
\vspace{1cm}
\end{tabular}
\end{table}}

Let $\alpha_i,h_i,\La_i$ ($i \in I$) be the simple roots, simple
coroots, and fundamental weights of $\geh$. Let $\delta$ and $c$
denote the generator of imaginary roots and the canonical central
element, respectively. Recall that $\delta=\sum_{i\in
I}a_i\alpha_i$ and $c = \sum_{i \in I}a^\vee_i h_i$, where the Kac
labels $a_i$ are the unique set of relatively prime positive
integers giving the linear dependency of the columns of the Cartan
matrix $A$ (that is, $A (a_0,\dotsc,a_n)^t = 0$). Explicitly,
\begin{equation}
\delta =
\begin{cases}
\alpha_0+\cdots + \alpha_n & \text{if $\geh=A^{(1)}_n$} \\
\alpha_0+\alpha_1+2\alpha_2+\cd+2\alpha_n & \text{if $\geh=B^{(1)}_{n}$} \\
\alpha_0+2\alpha_1+\cd+2\alpha_{n-1}+\alpha_n & \text{if $\geh=C^{(1)}_{n}$} \\
\alpha_0+\alpha_1+2\alpha_2+\cd+2\alpha_{n-2}
+\alpha_{n-1}+\alpha_n
& \text{if $\geh=D^{(1)}_{n}$} \\
2\alpha_0+2\alpha_1+\cd+2\alpha_{n-1}+\alpha_n
& \text{if $\geh=A^{(2)}_{2n}$} \\
\alpha_0+2\alpha_1+\cd+2\alpha_{n-1}+2\alpha_n
& \text{if $\geh=A^{(2)\dagger}_{2n}$} \\
\alpha_0+\alpha_1+2\alpha_2+\cd+2\alpha_{n-1}+\alpha_n & \text{if
$\geh=A^{(2)}_{2n-1}$} \\
\alpha_0+\alpha_1+\cd+\alpha_{n-1}+\alpha_n &\text{if
$\geh=D^{(2)}_{n+1}$}.
\end{cases}
\end{equation}
The dual Kac label $a^\vee_i$ is the label $a_i$ for the affine
Dynkin diagram obtained by ``reversing the arrows" of the Dynkin
diagram of $\geh$, or equivalently, the coefficients giving the
linear dependency of the rows of the Cartan matrix $A$. Note that
$a_0^\vee=2$ for $\geh=A^{(2)\dagger}_{2n}$ and $a_0^\vee=1$
otherwise.

Let $(\cdot|\cdot)$ be the normalized invariant form on $P$
\cite{Kac}. It satisfies
\begin{equation}
  (\alpha_i|\alpha_j) = \dfrac{a_i^\vee}{a_i} A_{ij}
\end{equation}
for $i,j\in I$. In particular
\begin{equation} \label{eq:form norm}
  (\alpha_a|\alpha_a)=\dfrac{2r}{a_0^\vee}
\end{equation}
if $\alpha_a$ is a long root.

For $i \in I$ let
\begin{equation}\label{eq:ttdef}
t_i = \max(\frac{a_i}{a^\vee_i},a^\vee_0), \qquad t^\vee_i =
\max(\frac{a^\vee_i}{a_i},a_0).
\end{equation}
The values $t_i$ are given in Table \ref{tab:Dynkin}. We shall
only use $t^\vee_i$ and $t_i$ for $i \in \I=I\backslash\{0\}$. For
$a\in\I$ we have
\begin{equation*}
t^\vee_a = 1 \,\,\text{ if $r = 1$,} \qquad t_a = a_0^\vee
\,\,\text{ if $r
> 1$}.
\end{equation*}

We consider two finite-dimensional subalgebras of $\geh$: $\gehb$,
whose Dynkin diagram is obtained from that of $\geh$ by removing
the $0$ vertex, and $\gf$, the subalgebra of $X_N$ fixed by the
automorphism $\sigma$ given in \cite[Section 8.3]{Kac}.
\begin{table}[ht]
\caption{}\label{tab:geh-bar}
\begin{center}
\begin{tabular}{c|cccccc}
$\geh$ & $X^{(1)}_N$ & $A^{(2)}_{2n}$ & $A^{(2)\dagger}_{2n}$ &
$A^{(2)}_{2n-1}$ &
$D^{(2)}_{n+1}$ \\
\hline
$\gehb$ & $X_N$ & $C_n$ & $B_n$ & $C_n$ & $B_n$ \\
$\gf$ &         $X_N$ & $B_n$ & $B_n$ & $C_n$ & $B_n$
\end{tabular}
\end{center}
\end{table}

Let $\gehb$ (resp. $\gf$) have weight lattice $\Pfin$ (resp.
$\Pt$), with simple roots and fundamental weights
$\alpha_a,\Lab_a$ (resp. $\alt_a,\Lt_a$) for $a\in\I$. Note that
$\gehb=\gf$ for $\geh\not=A^{(2)}_{2n}$. For $\geh=A^{(2)}_{2n}$,
$\gehb=C_n$ and $\gf=B_n$.

$\Pt$ is endowed with the bilinear form $(\cdot | \cdot )'$,
normalized by
\begin{equation} \label{eq:tilde form norm}
  (\alt_a|\alt_a )' = 2r/a_0^\vee \qquad\text{if
$\alt_a$ is a long root of $\gf$.}
\end{equation}
For $A^{(2)}_2$, the unique simple root $\alt_1$ of $\gf=B_1$ is
considered to be short.

Note that $\alpha_a,\Lab_a$ and $(\cdot | \cdot )$ may be
identified with $\alt_a,\Lt_a$ and $(\cdot | \cdot )'$ if $\geh
\not= A^{(2)}_{2n}$.

Define the $\Z$-linear map $\iota : \Pfin \rightarrow \Pt$ by
\begin{equation}\label{eq:iota}
\iota(\Lab_a) = \epsilon_a \Lt_a \qquad \text{for $a\in\I$,}
\end{equation}
where $\epsilon_a$ is defined by
\begin{equation}\label{eq:epsdef}
\epsilon_a = \begin{cases}
2 & \text{if $\geh=A^{(2)}_{2n}$ and $a=n$}\\
1 & \text{otherwise.}
\end{cases}
\end{equation}
In particular $\iota(\alpha_a) = \epsilon_a \alt_a$ for $a\in\I$.
We have
\begin{equation} \label{eq:forms}
(\iota(\alpha_b) | \iota(\alpha_b))' = a_0 (\alpha_b |
\alpha_b)\qquad\text{for all $b\in \I$.}
\end{equation}
If $\geh=A^{(2)}_{2n}$ both sides of \eqref{eq:forms} are equal to
$8$ if $b=n$ and $4$ otherwise. Especially for $\geh=A^{(2)}_2$
($n=1$), we have $(\tilde{\alpha}_1 | \tilde{\alpha}_1)' =2$ and
$(\alpha_1 | \alpha_1)=4$.
In the rest of the paper we shall write $(\cdot \vert\cdot)$ in
place of $(\cdot|\cdot)'$.

\subsection{Simple subalgebras}
\label{subsec:classical} For later use, specific realizations are
given for the simple roots and fundamental weights of the simple
Lie algebras of types $B_n$, $C_n$, and $D_n$, which appear as the
subalgebras $\gehb$ and $\gf$ of $\geh$. In each case the sublattice of
$\Pfin$ given by the weights appearing in tensor products of the
vector representation, is identified with $\Z^n$. Let
$\{\epsilon_i\mid 1\le i\le n\}$ be the standard basis of $\Z^n$.

\subsubsection*{The simple Lie algebra $B_n$}

\begin{align}
\label{eq:Broots}
\alpha_a&=\epsilon_a-\epsilon_{a+1} && \text{for $1\le a<n$}\\
\notag \alpha_n&=\epsilon_n &&\\
\notag
\Lab_a&=\epsilon_1+\cdots+\epsilon_a &&\text{for $1\le a<n$}\\
\notag \Lab_n&=\frac{1}{2}(\epsilon_1+\cdots+\epsilon_n).&&
\end{align}
$\la\in\Z^n$ is $B_n$-dominant if and only if
\begin{align} \label{eq:Bdom}
\la_a-\la_{a+1}&\ge0 &&\text{for $1\le a<n$} \\
\notag \la_n&\ge0.&
\end{align}

\subsubsection*{The simple Lie algebra $C_n$}

\begin{align} \label{eq:Croots}
\alpha_a&=\epsilon_a-\epsilon_{a+1} && \text{for $1\le a<n$}\\
\notag \alpha_n&=2\epsilon_n &&\\
\notag \Lab_a&=\epsilon_1+\cdots+\epsilon_a && \text{for $1\le
a\le n$.}
\end{align}
$\la\in\Z^n$ is $C_n$-dominant if and only if it is $B_n$-dominant
\eqref{eq:Bdom}.

\subsubsection*{The simple Lie algebra $D_n$}

\begin{align} \label{eq:Droots}
\alpha_a&=\epsilon_a-\epsilon_{a+1} &&\text{for $1\le a<n$}\\
\notag \alpha_n&=\epsilon_{n-1}+\epsilon_n &&\\
\notag
\Lab_a&=\epsilon_1+\cdots+\epsilon_a && \text{for $1\le a\le n-2$}\\
\notag \Lab_{n-1}&=\frac{1}{2}(\epsilon_1+\cdots+\epsilon_{n-1}-\epsilon_n)&&\\
\notag
\Lab_n&=\frac{1}{2}(\epsilon_1+\cdots+\epsilon_{n-1}+\epsilon_n)&&
\end{align}
$\la\in\Z^n$ is $D_n$-dominant if and only if
\begin{align} \label{eq:dominant D}
\la_a-\la_{a+1}&\ge 0 &&\text{for $1\le a<n$}\\
\notag \la_{n-1}+\la_n&\ge 0.&&
\end{align}

\subsection{Crystals}
Let $\geh'$ be the derived subalgebra of $\geh$. Denote the
corresponding quantized universal enveloping algebras of
$\geh\supset \geh'\supset \gehb$ by $U_q(\geh)\supset
U'_q(\geh)\supset U_q(\gehb)$.

In \cite{HKOTY} it is conjectured that there is a family of
finite-dimensional irreducible $U'_q(\geh)$-modules
$\{W^{(a)}_i\mid a\in\I,i\in \Z_{>0}\}$ which, unlike most
finite-dimensional $U'_q(\geh)$-modules, have crystal bases
$B^{a,i}$. This family is conjecturally characterized in several
different ways:
\begin{enumerate}
\item Its characters form the unique solutions of a system
of quadratic relations (the $Q$-system) \cite{KR1}.
\item Every crystal graph of an irreducible integrable finite-dimensional
$U'_q(\geh)$-module, is a tensor product of the $B^{a,i}$.
\item For $\la\in P$ let $V(\la)$ be the extremal weight module
defined in \cite[Section 3]{Ka} and $B(\la)$ its crystal base,
with unique vector $u_\la\in B(\la)$ of weight $\la$. Then the
affinization of $B^{a,i}$ (in the sense of \cite{KMN1}) is
isomorphic to the connected component of $u_\la$ in $B(\la)$, for
the weight $\la=i\Lab_a$ (except when $\geh=A^{(2)\dagger}_{2n}$
and $a=n$, in which case $\la=2i\Lab_a$).
\end{enumerate}

In light of point (2) above, we consider the category of crystal
graphs given by tensor products of the crystals $B^{a,i}$.

We introduce notation for tensor products of $B^{a,i}$. Let
$\mu=(L_i^{(a)})_{a\in\I,i\in\Z_{>0}}$ be a matrix of nonnegative
integers, almost all zero. Define
\begin{equation} \label{eq:path space}
 B^{(\mu)} = \bigotimes_{(a,i)\in\I\times\Z_{>0}}
(B^{a,i})^{\otimes L_i^{(a)}}.
\end{equation}
In type $A^{(1)}_n$ this is the tensor product of modules, which,
when restricted to $A_n$, are irreducible modules indexed by
rectangular partitions. The set of classically restricted paths
(or classical highest weight vectors) in $B^{(\mu)}$ of weight
$\la\in\Pfin^+=\bigoplus_{i\in\I} \Z_{\ge0} \Lab_i$ is by
definition
\begin{equation} \label{eq:cl hwv}
\Path(\la,\mu)=\{b\in B^{(\mu)}\mid \text{$\wt(b)=\la$ and $\et{i}b$
undefined for all $i\in \I$} \}.
\end{equation}
Here $\et{i}$ is given by the crystal graph. For $b,b'\in B^{a,i}$
we have $b'=\et{i}(b)$ if there is an arrow $b'\stackrel{i}{\longrightarrow}b$
in the crystal graph; if no such arrow exists then $\et{i}(b)$ is undefined.
Similarly, $b'=\ft{i}(b)$ if there is an arrow $b\stackrel{i}{\longrightarrow}b'$
in the crystal graph; if no such arrow exists then $\ft{i}(b)$ is undefined.
If $B_1$ and $B_2$ are crystals, then for $b_1\otimes b_2\in B_1\otimes B_2$
the action of $\et{i}$ is defined as
\begin{equation*}
\et{i}(b_1\otimes b_2)=\begin{cases}
\et{i}b_1 \otimes b_2 &\text{if $\varepsilon_i(b_1)>\varphi_i(b_2)$,}\\
b_1\otimes \et{i} b_2 &\text{else,}
\end{cases}
\end{equation*}
where $\varepsilon_i(b)=\max\{k\mid \et{i}^k \;\text{is defined}\}$ and
$\varphi_i(b)=\max\{k\mid \ft{i}^k \;\text{is defined}\}$.

\begin{assumption} \label{as:singlebox} In this paper we shall
restrict our attention to the case $B^{(\mu)}=B^{\otimes L}$ where
$B=B^{1,1}$. We shall write $B^{(\mut)}=B^{\otimes (L-1)}$.
\end{assumption}

The crystal graphs $B^{1,1}$ are listed in Table
\ref{tab:crystals}.
\begin{table}
\begin{tabular}{|c|l|}
\hline
%
$A_n^{(1)}$ & \raisebox{-0.7cm}{\scalebox{0.7}{
\begin{picture}(250,62)(-10,-12)
\BText(0,0){1} \LongArrow(10,0)(40,0) \BText(50,0){2}
\LongArrow(60,0)(90,0) \BText(100,0){3} \LongArrow(110,0)(140,0)
\Text(160,0)[]{$\cdots$} \LongArrow(175,0)(205,0)
\BText(220,0){n+1} \LongArrowArc(110,-181)(216,62,118)
\PText(25,2)(0)[b]{1} \PText(75,2)(0)[b]{2} \PText(125,2)(0)[b]{3}
\PText(190,2)(0)[b]{n} \PText(110,38)(0)[b]{0}
\end{picture}
}}
\\ \hline
%
$B_n^{(1)}$ & \raisebox{-1.3cm}{\scalebox{0.7}{
\begin{picture}(350,100)(-10,-50)
\BText(0,0){1} \LongArrow(10,0)(30,0) \BText(40,0){2}
\LongArrow(50,0)(70,0) \Text(85,0)[]{$\cdots$}
\LongArrow(95,0)(115,0) \BText(125,0){n} \LongArrow(135,0)(155,0)
\BText(165,0){0} \LongArrow(175,0)(195,0) \BBoxc(205,0)(13,13)
\Text(205,0)[]{\footnotesize$\overline{\mbox{n}}$}
\LongArrow(215,0)(235,0) \Text(250,0)[]{$\cdots$}
\LongArrow(260,0)(280,0) \BBoxc(290,0)(13,13)
\Text(290,0)[]{\footnotesize$\overline{\mbox{2}}$}
\LongArrow(300,0)(320,0) \BBoxc(330,0)(13,13)
\Text(330,0)[]{\footnotesize$\overline{\mbox{1}}$}
\LongArrowArc(185,-330)(365,68,112)
\LongArrowArcn(145,330)(365,-68,-112) \PText(20,2)(0)[b]{1}
\PText(60,2)(0)[b]{2} \PText(105,2)(0)[b]{n-1}
\PText(145,2)(0)[b]{n} \PText(185,2)(0)[b]{n}
\PText(225,2)(0)[b]{n-1} \PText(270,2)(0)[b]{2}
\PText(310,2)(0)[b]{1} \PText(185,38)(0)[b]{0}
\PText(145,-35)(0)[t]{0}
\end{picture}
}}
\\ \hline
%
%
$C_n^{(1)}$ & \raisebox{-0.7cm}{\scalebox{0.7}{
\begin{picture}(350,62)(-10,-12)
\BText(0,0){1} \LongArrow(10,0)(30,0) \BText(40,0){2}
\LongArrow(50,0)(70,0) \Text(85,0)[]{$\cdots$}
\LongArrow(95,0)(115,0) \BText(125,0){n} \LongArrow(135,0)(155,0)
\BBoxc(165,0)(13,13)
\Text(165,0)[]{\footnotesize$\overline{\mbox{n}}$}
\LongArrow(175,0)(195,0) \Text(210,0)[]{$\cdots$}
\LongArrow(220,0)(240,0) \BBoxc(250,0)(13,13)
\Text(250,0)[]{\footnotesize$\overline{\mbox{2}}$}
\LongArrow(260,0)(280,0) \BBoxc(290,0)(13,13)
\Text(290,0)[]{\footnotesize$\overline{\mbox{1}}$}
\LongArrowArc(145,-317)(352,67.5,112.5) \PText(20,2)(0)[b]{1}
\PText(60,2)(0)[b]{2} \PText(105,2)(0)[b]{n-1}
\PText(145,2)(0)[b]{n} \PText(185,2)(0)[b]{n-1}
\PText(230,2)(0)[b]{2} \PText(270,2)(0)[b]{1}
\PText(145,38)(0)[b]{0}
\end{picture}
}}
\\ \hline
%
$D_n^{(1)}$ & \raisebox{-1.3cm}{\scalebox{0.7}{
\begin{picture}(365,100)(-10,-50)
\BText(0,0){1} \LongArrow(10,0)(30,0) \BText(40,0){2}
\LongArrow(50,0)(70,0) \Text(85,0)[]{$\cdots$}
\LongArrow(95,0)(115,0) \BText(130,0){n-1}
\LongArrow(143,2)(160,14) \LongArrow(143,-2)(160,-14)
\BText(170,15){n} \BBoxc(170,-15)(13,13)
\Text(170,-15)[]{\footnotesize$\overline{\mbox{n}}$}
\LongArrow(180,14)(197,2) \LongArrow(180,-14)(197,-2)
\BBoxc(215,0)(25,13)
\Text(215,0)[]{\footnotesize$\overline{\mbox{n-1}}$}
\LongArrow(230,0)(250,0) \Text(265,0)[]{$\cdots$}
\LongArrow(275,0)(295,0) \BBoxc(305,0)(13,13)
\Text(305,0)[]{\footnotesize$\overline{\mbox{2}}$}
\LongArrow(315,0)(335,0) \BBoxc(345,0)(13,13)
\Text(345,0)[]{\footnotesize$\overline{\mbox{1}}$}
\LongArrowArc(192.5,-367)(402,69,111)
\LongArrowArcn(152.5,367)(402,-69,-111) \PText(20,2)(0)[b]{1}
\PText(60,2)(0)[b]{2} \PText(105,2)(0)[b]{n-2}
\PText(152,13)(0)[br]{n-1} \PText(152,-9)(0)[tr]{n}
\PText(188,13)(0)[bl]{n} \PText(188,-9)(0)[tl]{n-1}
\PText(240,2)(0)[b]{n-2} \PText(285,2)(0)[b]{2}
\PText(325,2)(0)[b]{1} \PText(192.5,38)(0)[b]{0}
\PText(152.5,-35)(0)[t]{0}
\end{picture}
}}
\\ \hline
%
$A_{2n}^{(2)}$ & \raisebox{-0.7cm}{\scalebox{0.7}{
\begin{picture}(350,62)(-10,-12)
\BText(0,0){1} \LongArrow(10,0)(30,0) \BText(40,0){2}
\LongArrow(50,0)(70,0) \Text(85,0)[]{$\cdots$}
\LongArrow(95,0)(115,0) \BText(125,0){n} \LongArrow(135,0)(155,0)
\BBoxc(165,0)(13,13)
\Text(165,0)[]{\footnotesize$\overline{\mbox{n}}$}
\LongArrow(175,0)(195,0) \Text(210,0)[]{$\cdots$}
\LongArrow(220,0)(240,0) \BBoxc(250,0)(13,13)
\Text(250,0)[]{\footnotesize$\overline{\mbox{2}}$}
\LongArrow(260,0)(280,0) \BBoxc(290,0)(13,13)
\Text(290,0)[]{\footnotesize$\overline{\mbox{1}}$}
\Text(145,35)[]{\footnotesize$\emptyset$}
\LongArrowArc(145,-317)(352,67.5,89)
\LongArrowArc(145,-317)(352,91,112.5) \PText(20,2)(0)[b]{1}
\PText(60,2)(0)[b]{2} \PText(105,2)(0)[b]{n-1}
\PText(145,2)(0)[b]{n} \PText(185,2)(0)[b]{n-1}
\PText(230,2)(0)[b]{2} \PText(270,2)(0)[b]{1}
\PText(220,30)(0)[b]{0} \PText(70,30)(0)[b]{0}
\end{picture}
}}
\\ \hline
%
%
$A_{2n}^{(2)\dagger}$ & \raisebox{-0.7cm}{\scalebox{0.7}{
\begin{picture}(350,62)(-10,-12)
\BText(0,0){1} \LongArrow(10,0)(30,0) \BText(40,0){2}
\LongArrow(50,0)(70,0) \Text(85,0)[]{$\cdots$}
\LongArrow(95,0)(115,0) \BText(125,0){n} \LongArrow(135,0)(155,0)
\BText(165,0){0} \LongArrow(175,0)(195,0) \BBoxc(205,0)(13,13)
\Text(205,0)[]{\footnotesize$\overline{\mbox{n}}$}
\LongArrow(215,0)(235,0) \Text(250,0)[]{$\cdots$}
\LongArrow(260,0)(280,0) \BBoxc(290,0)(13,13)
\Text(290,0)[]{\footnotesize$\overline{\mbox{2}}$}
\LongArrow(300,0)(320,0) \BBoxc(330,0)(13,13)
\Text(330,0)[]{\footnotesize$\overline{\mbox{1}}$}
\PText(20,2)(0)[b]{1} \PText(60,2)(0)[b]{2}
\PText(105,2)(0)[b]{n-1} \PText(145,2)(0)[b]{n}
\PText(185,2)(0)[b]{n} \PText(225,2)(0)[b]{n-1}
\PText(270,2)(0)[b]{2} \PText(310,2)(0)[b]{1}
\PText(165,38)(0)[b]{0}
\LongArrowArc(165,-433)(468,71,109)
\end{picture}
}}
\\ \hline
%
$A_{2n-1}^{(2)}$ & \raisebox{-1.3cm}{\scalebox{0.7}{
\begin{picture}(310,100)(-10,-50)
\BText(0,0){1} \LongArrow(10,0)(30,0) \BText(40,0){2}
\LongArrow(50,0)(70,0) \Text(85,0)[]{$\cdots$}
\LongArrow(95,0)(115,0) \BText(125,0){n} \LongArrow(135,0)(155,0)
\BBoxc(165,0)(13,13)
\Text(165,0)[]{\footnotesize$\overline{\mbox{n}}$}
\LongArrow(175,0)(195,0) \Text(210,0)[]{$\cdots$}
\LongArrow(220,0)(240,0) \BBoxc(250,0)(13,13)
\Text(250,0)[]{\footnotesize$\overline{\mbox{2}}$}
\LongArrow(260,0)(280,0) \BBoxc(290,0)(13,13)
\Text(290,0)[]{\footnotesize$\overline{\mbox{1}}$}
\LongArrowArc(165,-240)(275,65,115)
\LongArrowArcn(125,240)(275,-65,-115) \PText(20,2)(0)[b]{1}
\PText(60,2)(0)[b]{2} \PText(105,2)(0)[b]{n-1}
\PText(145,2)(0)[b]{n} \PText(185,2)(0)[b]{n-1}
\PText(230,2)(0)[b]{2} \PText(270,2)(0)[b]{1}
\PText(165,38)(0)[b]{0} \PText(125,-35)(0)[t]{0}
\end{picture}
}}
\\ \hline
%
$D_{n+1}^{(2)}$ & \raisebox{-0.7cm}{\scalebox{0.7}{
\begin{picture}(350,62)(-10,-12)
\BText(0,0){1} \LongArrow(10,0)(30,0) \BText(40,0){2}
\LongArrow(50,0)(70,0) \Text(85,0)[]{$\cdots$}
\LongArrow(95,0)(115,0) \BText(125,0){n} \LongArrow(135,0)(155,0)
\BText(165,0){0} \LongArrow(175,0)(195,0) \BBoxc(205,0)(13,13)
\Text(205,0)[]{\footnotesize$\overline{\mbox{n}}$}
\LongArrow(215,0)(235,0) \Text(250,0)[]{$\cdots$}
\LongArrow(260,0)(280,0) \BBoxc(290,0)(13,13)
\Text(290,0)[]{\footnotesize$\overline{\mbox{2}}$}
\LongArrow(300,0)(320,0) \BBoxc(330,0)(13,13)
\Text(330,0)[]{\footnotesize$\overline{\mbox{1}}$}
\Text(165,35)[]{\footnotesize$\emptyset$}
\LongArrowArc(165,-433)(468,71,89)
\LongArrowArc(165,-433)(468,91,109) \PText(20,2)(0)[b]{1}
\PText(60,2)(0)[b]{2} \PText(105,2)(0)[b]{n-1}
\PText(145,2)(0)[b]{n} \PText(185,2)(0)[b]{n}
\PText(225,2)(0)[b]{n-1} \PText(270,2)(0)[b]{2}
\PText(310,2)(0)[b]{1} \PText(248,30)(0)[b]{0}
\PText(82,30)(0)[b]{0}
\end{picture}
}}
\\ \hline
\end{tabular}\vspace{4mm}
\caption{\label{tab:crystals}Crystals $B^{1,1}$}
\end{table}

In each case (other than $A^{(1)}_n$) the elements of $B=B^{1,1}$
consist of $\{k,\overline{k}\mid 1\le k\le n\}$ and possibly
elements $0$ and $\phi$.

\begin{remark} \label{rem:highest weight}
By glancing at Table \ref{tab:crystals}, one may check that the
following are equivalent for $b=b_L\otimes
b_{L-1}\otimes\dotsm\otimes b_1\in B^{\otimes L}$ and
$\la\in\Pfin^+$.
\begin{enumerate}
\item $b$ is a classically restricted path of weight
$\la\in\Pfin^+$.
\item $\la-\wt(b_L)\in\Pfin^+$, $b_{L-1}\otimes\dotsm\otimes b_1$
is a classically restricted path of weight $\la-\wt(b_L)$, and if
$b_L=0\in B$ then $\la_n>0$ (where $\la$ is viewed as an element
of $\Z^n$).
\end{enumerate}
\end{remark}

The weight function $\wt:B\to \Z^n$ is given by
\begin{alignat*}{2}
  \wt(k) &= \epsilon_k \qquad&\text{for $1\le k\le n$} \\
  \wt(\overline{k}) &= -\epsilon_k\qquad&\text{for $1\le k\le n$}
  \\
  \wt(0) &= \wt(\phi) = 0.
\end{alignat*}
The weight function $\wt:B^{\otimes L}\to\Z^n$ is defined by
$\wt(b_L\otimes\dotsm\otimes b_1) =\sum_{j=1}^L \wt(b_j)$. So if
$\la=\wt(p)$ where $p\in B^{\otimes L}$, then $\la_k$ is the
multiplicity of $k$ in $p$ minus the multiplicity of
$\overline{k}$ in $p$.

\subsection{One-dimensional sums}
\label{subsec:1dsum} The energy function
$D:B^{(\mu)}\rightarrow\Z$ gives the grading on $B^{(\mu)}$. In
the case $B^{(\mu)}=B^{\otimes L}$ it takes a simple form. Due to
the existence of the universal $R$-matrix and the fact that
$W^{(1)}_1$ is irreducible, by \cite{KMN1} there is a unique (up
to global additive constant) function $H:B^{1,1}\otimes
B^{1,1}\rightarrow\Z$ called the local energy function, such that
\begin{equation} \label{eq:loc en}
  H(\et{i}(b\otimes b')) = H(b\otimes b')+
  \begin{cases}
  -1 & \text{if $i=0$ and $\et{0}(b\otimes b')=b\otimes\et{0}b'$}
  \\
  1 & \text{if $i=0$ and $\et{0}(b\otimes b')=\et{0}b\otimes b'$}
  \\
  0 & \text{otherwise.}
  \end{cases}
\end{equation}
Let $b^\natural\in B^{1,1}$ be the unique element such that
$\vphi(b^\natural)=\La_0$. We normalize $H$ by the condition
\begin{equation} \label{eq:H norm}
  H(1 \otimes 1) = 0.
\end{equation}
Then
\begin{equation} \label{eq:energy}
\begin{split}
  E(b_L\otimes\dotsm \otimes b_1) &= L\,\, H(b_1\otimes b^\natural)+
  \sum_{j=1}^{L-1} (L-j) \,\,H(b_{j+1}\otimes b_j),\\
  D(b_L\otimes\dotsm \otimes b_1) &= E(b_L\otimes\dotsm \otimes b_1)
   - E(1\otimes\cdots \otimes 1).
\end{split}
\end{equation}
Define the one-dimensional sum $X(\la,\mu;q)\in\Z[q,q^{-1}]$ by
\begin{equation} \label{eq:onedim}
  X(\la,\mu;q) = \sum_{b\in \Path(\la,\mu)} q^{D(b)}.
\end{equation}
Since $B^{(\mu)}$ is completely reducible as a
$U_q(\gehb)$-crystal, one has
\begin{equation*}
  \sum_{b\in B^{(\mu)}} e^{\wt(b)} q^{D(b)} =
  \sum_{\la\in\Pfin^+} \chi^\la X(\la,\mu;q)
\end{equation*}
where $\chi^\la$ is the character of the irreducible
$U_q(\gehb)$-module of highest weight $\la$. It can be shown that
$X(\la,\mu;q)\in\Z_{\ge0}[q^{-1}]$. For convenience we define
\begin{equation} \label{eq:co}
    \Hb = - H, \qquad \Db = - D,\qquad
    \Xb(\la,\mu;q)=X(\la,\mu;q^{-1}).
\end{equation}

\section{Rigged configurations and the bijection}
\label{sec:bijection}
\subsection{The fermionic formula, $\geh\not=A^{(2)\dagger}_{2n}$}
\label{subsec:ferm} This subsection reviews definitions of
\cite{HKOTT,HKOTY}. Let $\geh$ be a Kac-Moody algebra of
nonexceptional affine type that is not of the form
$A^{(2)\dagger}_{2n}$. Fix $\la\in\Pfin^+$ and a matrix
$\mu=(L_i^{(a)})$ of nonnegative integers as in subsection
\ref{subsec:1dsum}.

Let $\nu=(m_i^{(a)})$ be another such matrix. Say that $\nu$ is a
$\la$-configuration if
\begin{equation}
\label{eq:config} \sum_{\substack{a\in\I \\ i\in\Z_{>0}}} i\,
m_i^{(a)} \alt_a = \iota\left( \sum_{\substack{a\in\I \\
i\in\Z_{>0}}} i \,L_i^{(a)} \Lab_a - \la \right).
\end{equation}
Say that a configuration $\nu$ is $\mu$-admissible if
\begin{equation} \label{eq:ppos}
  p_i^{(a)} \ge 0\qquad\text{for all $a\in\I$ and
  $i\in\Z_{>0}$,}
\end{equation}
where
\begin{equation} \label{eq:p}
p_i^{(a)} = \sum_{k\in\Z_{>0}} \left( L_k^{(a)} \min(i,k) -
\dfrac{1}{t_a^\vee} \sum_{b\in\I} (\alt_a|\alt_b)\min(t_b i,t_a
k)\, m_k^{(b)}\right).
\end{equation}
Write $C(\la,\mu)$ for the set of $\mu$-admissible
$\la$-configurations. Define
\begin{equation} \label{eq:cc}
cc(\nu) = \dfrac{1}{2} \sum_{a,b\in\I} \sum_{j,k\in\Z_{>0}} (\alt_a|\alt_b)
\min(t_b j, t_a k) m_j^{(a)} m_k^{(b)}.
\end{equation}
The fermionic formula is defined by
\begin{equation}\label{fermi}
\Mb(\la,\mu;q) = \sum_{\nu\in C(\la,\mu)} q^{cc(\nu)}
\prod_{a\in\I} \prod_{i\in\Z_{>0}}
\qbin{p_i^{(a)}+m_i^{(a)}}{m_i^{(a)}}_{q^{t^\vee_a}}.
\end{equation}

The $X=M$ conjecture of \cite{HKOTT,HKOTY} states that
\begin{equation}\label{eq:X=M}
  \Xb(\la,\mu;q)=\Mb(\la,\mu;q).
\end{equation}

\subsection{Rigged configurations, $\geh\not=A^{(2)\dagger}_{2n}$}
The fermionic formula $\Mb(\la,\mu)$ can be interpreted using
combinatorial objects called rigged configurations. These objects
are a direct combinatorialization of the fermionic formula
$\Mb(\la,\mu;q)$. Our goal is to prove \eqref{eq:X=M} under
Assumption \ref{as:singlebox} by defining a statistic-preserving
bijection from rigged configurations to paths. For this purpose it
is convenient to use an indexing slightly differing from that used
above.

For $a\in\I$, define
\begin{equation}
  \f_a = \begin{cases}
  2 & \text{if $a=n$ and $\geh=C_n^{(1)}$} \\
  \frac{1}{2} & \text{if $a=n$ and $\geh=B_n^{(1)}$} \\
  1 & \text{otherwise.}
  \end{cases}
\end{equation}
$\f_a$ is half the square length of $\alpha_a$ for untwisted
affine types and is equal to $1$ for twisted types.

A quasipartition $\la$ of type $a\in\I$ is a finite multiset taken
from the set $\f_a \Z_{>0}$. Denote by $m_i(\la)$ the number of
times $i\in\f_a\Z_{>0}$ occurs in $\la$. The diagram of such a
quasipartition has, for each $i\in\f_a\Z_{>0}$, $m_i(\la)$ rows
consisting of $i$ boxes, where each box has width $\f_a$. Set
\begin{equation} \label{eq:index pairs}
  \HHH=\{(a,i)\mid a\in\I, i\in\f_a\Z_{>0}\}.
\end{equation}
Denote by $(\nt,\Jt)$ a pair where $\nt=\{\nu^{(a)}\}_{a\in \I}$
is a sequence of quasipartitions with $\nu^{(a)}$ of type $a$ and
$\Jt=\{J^{(a,i)}\}_{(a,i)\in\HHH}$ is a double sequence of
partitions. For $(a,i)\in\HHH$, define
\begin{equation}\label{eq:correspondence}
\begin{split}
P_i^{(a)}(\nt) &= p_{i/\f_a}^{(a)} \\
m_i^{(a)}(\nt) &= m_{i/\f_a}^{(a)}=m_i(\nu^{(a)}).
\end{split}
\end{equation}
Then a rigged configuration is a pair $(\nt,\Jt)$ subject to the
restriction \eqref{eq:config} and the requirement that $J^{(a,i)}$
be a quasipartition contained in a $m_i^{(a)}(\nt) \times
P_i^{(a)}(\nt)$ rectangle. The set of rigged configurations for
fixed $\la$ and $\mu$ is denoted by $\RC(\la,\mu)$. Then
\eqref{fermi} is equivalent to
\begin{equation*}
F(\la,\mu)=\sum_{(\nt,\Jt)\in\RC(\la,\mu)} q^{cc(\nt,\Jt)}
\end{equation*}
where $cc(\nt,\Jt)=cc(\nu)+|\Jt|$
and $|\Jt|=\sum_{(a,i)\in \HHH} t_a^\vee |J^{(a,i)}|$
for $\nu$ corresponding to $\nt$ under \eqref{eq:correspondence}.

\subsection{$A^{(2)\dagger}_{2n}$ rigged configurations} In this
subsection let $\geh=A^{(2)\dagger}_{2n}$. As this case is not
considered in \cite{HKOTT} we shall only give the definition in
terms of rigged configurations, although it is easy to express the
result as a sum of a product of $q$-binomials (see \cite[Section 7.6]{OSS}).
The important feature
is that the riggings of odd-sized parts of $\nu^{(n)}$, must have
the form $x/2$ where $x$ is an odd integer. So let $\mu$ and $\la$
be as in subsection \ref{subsec:ferm}. Given a matrix
$\nu=(m^{(a)}_i)$, let $P_i^{(a)}(\nt)$ and $m_i^{(a)}(\nt)$ be
defined as before. Call $\nt$ $\mu$-admissible if
$P_i^{(a)}(\nt)\ge0$ for all $a\in\I$ and $i\in\Z_{>0}$, together
with the extra condition that
\begin{equation} \label{eq:a2doddvac}
  P_i^{(n)}(\nt) \ge 1 \qquad\text{if $i$ is odd and $m_i^{(n)}(\nt)>0$.}
\end{equation}
A rigging $\Jt$ consists of quasipartitions $J^{(a,i)}$ for
$a\in\I$ and $i\in\Z_{>0}$. For $a\not=n$ or $i$ even, $J^{(a,i)}$
is an ordinary partition satisfying the usual properties. For
$a=n$ and $i$ odd, $J^{(n,i)}$ is a quasipartition contained in a
rectangle with $P_i^{(n)}(\nt)$ columns and $m_i^{(n)}(\nt)$ rows,
but it has cells of width 1/2 and each part size must be of the
form $x/2$ for $x$ an odd integer. This defines the set
$RC(\la,\mu)$ for $\geh=A^{(2)\dagger}_{2n}$. Then $F(\la,\mu)$ is
defined as before where $|\Jt|$ is the sum of the areas of all the
quasipartitions $J^{a,i}$. This definition is compatible with the
virtual crystal realization which embeds paths (and rigged
configurations) of type $A^{(2)\dagger}_{2n}$ into those of type
$A^{(1)}_{2n-1}$ \cite{OSS}.

\subsection{The bijection from RCs to paths}
We now describe the general form of the bijection
$\Phi:\RC(\la,\mu)\to\Path(\la,\mu)$ under Assumption
\ref{as:singlebox}. Let $\mu=(L_i^{(a)})$ be such that
$B^{(\mu)}=B^{\otimes L}$, that is, $L_i^{(a)}=L \delta_{a1}
\delta_{i1}$. Let $\mut$ be such that $B^{(\mut)}=B^{\otimes
(L-1)}$.

Let $(\nt,\Jt)\in\RC(\la,\mu)$. We shall define a map
$\rk:\RC(\la,\mu)\to B$ which associates to $(\nt,\Jt)$ an element
of $B$ called its rank.

Denote by $\RC_b(\la,\mu)$ the elements of $\RC(\la,\mu)$ of rank
$b$. We shall define a bijection
$\delta:\RC_b(\la,\mu)\to\RC(\la-\wt(b),\mut)$. The disjoint union
of these bijections then defines a bijection
$\delta:\RC(\la,\mu)\to\bigcup_{b\in B} \RC(\la-\wt(b),\mut)$.

The bijection $\Phi$ is defined recursively as follows. For $b\in
B$ let $\Path_b(\la,\mu)$ be the set of paths in
$B^{(\mu)}=B^{\otimes L}$ that have $b$ as leftmost tensor factor.
For $L=0$ the bijection $\Phi$ sends the empty rigged
configuration (the only element of the set $\RC(\la,\mu)$) to the
empty path (the only element of $\Path(\la,\mu)$). Otherwise
assume that $\Phi$ has been defined for $B^{\otimes (L-1)}$ and
define it for $B^{\otimes L}$ by the commutative diagram
\begin{equation}
\begin{CD}
\RC_b(\la,\mu) @>{\Phi}>> \Path_b(\la,\mu) \\
@V{\delta}VV @VVV \\
\RC(\la-\wt(b),\mut) @>{\Phi}>> \Path(\la-\wt(b),\mut)
\end{CD}
\end{equation}
where the right hand vertical map removes the leftmost tensor
factor $b$. In short,
\begin{equation}
  \Phi(\nt,\Jt)=\rk(\nt,\Jt)\otimes \Phi(\delta(\nt,\Jt)).
\end{equation}

\begin{remark} \label{rem:phidef} For $\Phi$ to be well-defined,
by Remark \ref{rem:highest weight} it must be shown that if
$b=\rk(\nt,\Jt)$, then $\rho=\la-\wt(b)$ is dominant, and if $b=0$
then $\la_n>0$.
\end{remark}

We also require the bijection
$\Phit:\RC(\la,\mu)\to\Path(\la,\mu)$ given by $\Phit=\Phi\circ
\comp$ where $\comp:\RC(\la,\mu)\to\RC(\la,\mu)$ with
$\comp(\nt,\Jt)=(\nt,\Jtt)$ is the function which complements the
riggings, meaning that $\Jtt$ is obtained from $\Jt$ by
complementing all partitions $J^{(a,i)}$ in the $m_i^{(a)}\times
P_i^{(a)}(\nt)$ rectangle.

\begin{theorem}\label{thm:bij}
$\Phi:\RC(\la,\mu)\to\Path(\la,\mu)$ is a bijection such that
\begin{equation} \label{eq:cc=D}
cc(\nt,\Jt)=\Db(\Phit(\nt,\Jt))\qquad\text{for all
$(\nt,\Jt)\in\RC(\la,\mu)$.}
\end{equation}
\end{theorem}
For type $A_n^{(1)}$ a generalization of this theorem for
all $\mu$ was proven in \cite{KSS}.
For other types Theorem \ref{thm:bij} is proved in section
\ref{sec:proof}.

\section{The bijection for each root system}
In this section the maps $\rk$ and $\delta$ are defined in a
case-by-case manner. For each $\geh$, an explicit formula is given
for the vacancy numbers $P^{(a)}_i(\nt)$ (see
\eqref{eq:correspondence}), obtained by writing \eqref{eq:p} in
terms of the function $Q_i$ (see \eqref{eq:Qdef}) using the data
for the simple Lie algebras given in section
\ref{subsec:classical}. Then for $(\nt,\Jt)\in\RC(\la,\mu)$, an
algorithm is given which defines $b=\rk(\nt,\Jt)$, the new smaller
rigged configuration $(\ntt,\Jtt)=\delta(\nt,\Jt)$ such that
$(\ntt,\Jtt)\in\RC(\rho,\mut)$ (where $\rho=\la-\wt(b)$), and the
new vacancy numbers in terms of the old.

For a quasipartition $\tau$ with boxes of width $\f$ and
$i\in\f\N$, define
\begin{equation} \label{eq:Qdef}
Q_i(\tau)=\sum_j \min(\tau_j,i),
\end{equation}
the area of $\tau$ in the first $i$ quasicolumns.

The quasipartition $J^{(a,i)}$ is called \textit{singular} (with
respect to the configuration $\nt$) if it has a part of size
$P_i^{(a)}(\nt)$. If $A$ is a statement then $\chi(A)=1$ if $A$ is
true and $\chi(A)=0$ if $A$ is false. We also use the Kronecker delta
notation $\delta_{a,b}=\chi(a=b)$.

\subsection{Bijection algorithm for type $D_n^{(1)}$}
\subsection*{Vacancy numbers}
\begin{equation}\label{eq:vac D}
\begin{split}
P_i^{(a)}(\nt)&=Q_i(\nu^{(a-1)})-2Q_i(\nu^{(a)})+Q_i(\nu^{(a+1)})
 +L\delta_{a,1} \qquad\text{for $1\le a<n-2$}\\
P_i^{(n-2)}(\nt)&=Q_i(\nu^{(n-3)})-2Q_i(\nu^{(n-2)})+Q_i(\nu^{(n-1)})
 +Q_i(\nu^{(n)})\\
P_i^{(n-1)}(\nt)&=Q_i(\nu^{(n-2)})-2Q_i(\nu^{(n-1)})\\
P_i^{(n)}(\nt)&=Q_i(\nu^{(n-2)})-2Q_i(\nu^{(n)})
\end{split}
\end{equation}
\subsection*{Constraints}
\begin{equation}\label{eq:constraint D}
\begin{aligned}[2]
|\nu^{(a)}|&=L-\sum_{b=1}^a \la_b &\text{for $1\le a\le n-2$}\\
|\nu^{(n-1)}|&=\frac{1}{2}(L-\sum_{b=1}^{n-1} \la_b+\la_n)&\\
|\nu^{(n)}|&=\frac{1}{2}(L-\sum_{b=1}^n \la_b)&
\end{aligned}
\end{equation}
\subsection*{Algorithm $\delta$}
Set $\ell^{(0)}=0$ and repeat the following process for
$a=1,2,\ldots,n-2$ or until stopped. Find the minimal index $i\ge
\ell^{(a-1)}$ such that $J^{(a,i)}$ is singular. If no such $i$
exists, set $b=a$ and stop. Otherwise set $\ell^{(a)}=i$ and
continue with $a+1$.

If the process has not stopped at $a=n-2$ continue as follows.
Find the minimal indices $i,j\ge \ell^{(n-2)}$ such that
$J^{(n-1,i)}$ and $J^{(n,j)}$ are singular. If neither $i$ nor
$j$ exist, set $b=n-1$ and stop.
If $i$ exists, but not $j$, set $\ell^{(n-1)}=i$, $b=n$ and stop.
If $j$ exists, but not $i$, set $\ell^{(n)}=j$, $b=\overline{n}$
and stop. If both $i$ and $j$ exist, set $\ell^{(n-1)}=i$, $\ell^{(n)}=j$
and continue with $a=n-2$.

Now continue for $a=n-2,n-3,\ldots,1$ or until stopped.
Find the minimal index $i\ge \lb^{(a+1)}$ where $\lb^{(n-1)}
=\max(\ell^{(n-1)},\ell^{(n)})$ such that $J^{(a,i)}$ is singular
(if $i=\ell^{(a)}$ then there need to be two parts of size
$P_i^{(a)}(\nt)$ in $J^{(a,i)}$).
If no such $i$ exists, set $b=\overline{a+1}$ and stop.
If the process did not stop, set $b=\overline{1}$.

Set all yet undefined $\ell^{(a)}$ and $\lb^{(a)}$ to $\infty$.
\subsection*{New RC}
\begin{equation} \label{eq:ch m D1}
 m_i^{(a)}(\ntt)=m_i^{(a)}(\nt)+\begin{cases}
 1 & \text{if $i=\ell^{(a)}-1$}\\
 -1 & \text{if $i=\ell^{(a)}$}\\
 1 & \text{if $i=\lb^{(a)}-1$ and $1\le a\le n-2$}\\
 -1 & \text{if $i=\lb^{(a)}$ and $1\le a \le n-2$}\\
 0 & \text{otherwise} \end{cases}
\end{equation}
The partition $\tilde{J}^{(a,i)}$ is obtained from $J^{(a,i)}$ by removing
a part of size $P_i^{(a)}(\nt)$ for $i=\ell^{(a)}$ and $i=\lb^{(a)}$,
adding a part of size $P_i^{(a)}(\ntt)$ for $i=\ell^{(a)}-1$ and
$i=\lb^{(a)}-1$, and leaving it unchanged otherwise.

\subsection*{Change in vacancy numbers}
\begin{align}\label{eq:cv D}
P_i^{(a)}(\ntt)=P_i^{(a)}(\nt)-&\chi(\ell^{(a-1)}\le i)
 +2\chi(\ell^{(a)}\le i)-\chi(\ell^{(a+1)}\le i)\\
 -&\chi(\lb^{(a-1)}\le i)+2\chi(\lb^{(a)}\le i)-\chi(\lb^{(a+1)}\le i)\notag\\
\intertext{for $1\le a<n-2$}
P_i^{(n-2)}(\ntt)=P_i^{(n-2)}(\nt)-&\chi(\ell^{(n-3)}\le i)
 +2\chi(\ell^{(n-2)}\le i)-\chi(\ell^{(n-1)}\le i)\notag\\
 -&\chi(\lb^{(n-3)}\le i)+2\chi(\lb^{(n-2)}\le i)-\chi(\ell^{(n)}\le i)\notag\\
P_i^{(n-1)}(\ntt)=P_i^{(n-1)}(\nt)-&\chi(\ell^{(n-2)}\le i)
 -\chi(\lb^{(n-2)}\le i)+2\chi(\ell^{(n-1)}\le i)\notag\\
P_i^{(n)}(\ntt)=P_i^{(n)}(\nt)-&\chi(\ell^{(n-2)}\le i)
 -\chi(\lb^{(n-2)}\le i)+2\chi(\ell^{(n)}\le i).\notag
\end{align}

\subsection{Bijection algorithm for type $B_n^{(1)}$}

\subsection*{Vacancy numbers}
\begin{equation}\label{eq:vac B}
\begin{aligned}[2]
P_i^{(a)}(\nt)&=Q_i(\nu^{(a-1)})-2Q_i(\nu^{(a)})+Q_i(\nu^{(a+1)})
 +L\delta_{a,1}
 && \text{for $i\in\N$} \\
 &&&\text{$1\le a\le n-2$}\\
P_i^{(n-1)}(\nt)&=Q_i(\nu^{(n-2)})-2Q_i(\nu^{(n-1)})+2Q_i(\nu^{(n)})
 && \text{for $i\in \N$}\\
P_i^{(n)}(\nt)&=2Q_i(\nu^{(n-1)})-4Q_i(\nu^{(n)})
 && \text{for $i\in \frac{1}{2} \N$}
\end{aligned}
\end{equation}

\subsection*{Constraints}
\begin{equation}\label{eq:constraint B}
\begin{aligned}[2]
|\nu^{(a)}|&=L-\sum_{b=1}^a \la_b &\text{for $1\le a\le n-1$}\\
|\nu^{(n)}|&=\frac{1}{2}(L-\sum_{b=1}^n \la_b)&
\end{aligned}
\end{equation}

\subsection*{Algorithm $\delta$}
Call a partition \textit{quasi-singular} if it is not singular and has a
part of size $P_i^{(a)}(\nt)-1$.

Set $\ell^{(0)}=0$ and repeat the following process for
$a=1,2,\ldots,n-1$ or until stopped.
Find the minimal index $i\ge \ell^{(a-1)}$ such that $J^{(a,i)}$ is singular.
If no such $i$ exists, set $b=a$ and stop.
Otherwise set $\ell^{(a)}=i$ and continue.

If the process has not yet stopped, continue as follows.
For brevity let us denote by (S) and (Q) the following conditions:
\begin{enumerate}
\item[(S)] $i\ge \ell^{(n-1)}$ and $J^{(n,i)}$ is singular.
\item[(Q)] $i=\ell^{(n-1)}-\frac{1}{2}$ and $J^{(n,i)}$ is singular; or
           $i\ge \ell^{(n-1)}$ and $J^{(n,i)}$ is quasi-singular.
\end{enumerate}
Find the minimal index $i\ge \ell^{(n-1)}-\frac{1}{2}$ such that
(S) or (Q) holds (note that (S) and (Q) are mutually excluding).
If no such $i$ exists, set $b=n$ and stop.
If (S) holds set $\lb^{(n)}=i$ and $\ell^{(n)}=i-\frac{1}{2}$. Say that
case (S) holds.
If (Q) holds set $\ell^{(n)}=i$ and find the minimal index
$j>i$ such that (S) holds. If no such $j$ exists, set $b=0$ and stop.
Say that case (Q) holds.
Otherwise set $\lb^{(n)}=j$ and say that case (Q,S) holds.

If the process has not yet stopped continue in the following
fashion for $a=n-1,n-2,\ldots,1$ or until stopped. Find the
minimal index $i\ge \lb^{(a+1)}$ such that $J^{(a,i)}$ is singular
(if $\ell^{(a)}=i$ then $J^{(a,i)}$ actually needs to have two
parts of size $P_i^{(a)}(\nt)$). If no such $i$ exists, set
$b=\overline{a+1}$ and stop. Otherwise set $\lb^{(a)}=i$ and
continue. If the process did not stop for $a\ge 1$ set
$b=\overline{1}$.

Set all undefined $\ell^{(a)}$ and $\lb^{(a)}$ for $1\le a\le n$
to $\infty$.
\subsection*{New RC}
\begin{equation} \label{eq:ch m B1}
m_i^{(a)}(\ntt)=m_i^{(a)}(\nt)+\begin{cases}
 1 & \text{if $i=\ell^{(a)}-\f_a$}\\
 -1 & \text{if $i=\ell^{(a)}$}\\
 1 & \text{if $i=\lb^{(a)}-\f_a$}\\
 -1 & \text{if $i=\lb^{(a)}$}\\
 0 & \text{otherwise} \end{cases}
\end{equation}
Note that if two or more conditions hold, all of the changes should
be performed.

For $1\le a<n$ the partition $\tilde{J}^{(a,i)}$ is obtained from
$J^{(a,i)}$ by removing a part of size $P_i^{(a)}(\nt)$ for $i=\ell^{(a)}$
and $i=\lb^{(a)}$, adding a part of size $P_i^{(a)}(\ntt)$ for
$i=\ell^{(a)}-1$ and $i=\lb^{(a)}-1$ and leaving it unchanged
otherwise.
If case (S) occurred $\tilde{J}^{(n,i)}$ is obtained from $J^{(n,i)}$ by
removing a part of size $P_i^{(n)}(\nt)$ for $i=\lb^{(n)}$, adding a part
of size $P_i^{(n)}(\ntt)$ for $i=\lb^{(n)}-1$, and leaving it unchanged
otherwise. If case (Q) holds remove the largest part in $J^{(n,i)}$
for $i=\ell^{(n)}$ and add a part of size $P_i^{(n)}(\ntt)$ for
$i=\ell^{(n)}-\frac{1}{2}$. If case (Q,S) holds, then apply (S') for
$t=\ell^{(n)}$ and (Q') for $t=\lb^{(n)}$ where
\begin{enumerate}
\item[(S')] obtain $\tilde{J}^{(n,i)}$ from $J^{(n,i)}$
 by removing the largest part for $i=t$ and adding a part of
 size $P_i^{(n)}(\ntt)$ for $i=t-\frac{1}{2}$, leaving all
 other $J^{(n,i)}$ unchanged;
\item[(Q')] obtain $\tilde{J}^{(n,i)}$ from $J^{(n,i)}$
 by removing the largest part for $i=t$ and adding a part of
 size $P_i^{(n)}(\ntt)-1$ if $t<\lb^{(n-1)}$ and of size $P_i^{(n)}(\ntt)$ if
 $t=\lb^{(n-1)}$ for $i=t-\frac{1}{2}$, leaving all other $J^{(n,i)}$
 unchanged.
\end{enumerate}

\subsection*{Change in vacancy numbers}
\begin{align}\label{eq:cv B}
P_i^{(a)}(\ntt)=P_i^{(a)}(\nt)&-\chi(\ell^{(a-1)}\le i)
 +2\chi(\ell^{(a)}\le i)-\chi(\ell^{(a+1)}\le i)\\
 &-\chi(\lb^{(a-1)}\le i)
 +2\chi(\lb^{(a)}\le i)-\chi(\lb^{(a+1)}\le i)\notag \\
\intertext{for $1\le a\le n-1$ and}
P_i^{(n)}(\ntt)=P_i^{(n)}(\nt)&-\chi(\ell^{(n-1)}-\frac{1}{2}\le i)
 -\chi(\ell^{(n-1)}\le i)+2\chi(\ell^{(n)}\le i)\notag\\
 &-\chi(\lb^{(n-1)}-\frac{1}{2}\le i)
 -\chi(\lb^{(n-1)}\le i)+2\chi(\lb^{(n)}\le i).\notag
\end{align}

\subsection{Bijection algorithm for type $C_n^{(1)}$}

\subsection*{Vacancy numbers}
\begin{equation}\label{eq:vac C}
\begin{aligned}[3]
P_i^{(a)}(\nt)&=Q_i(\nu^{(a-1)})-2Q_i(\nu^{(a)})+Q_i(\nu^{(a+1)})
 +L\delta_{a,1} && \text{for $i\in\Z_{\ge 0}$}\\
&&& 1\le a<n\\
P_i^{(n)}(\nt)&=Q_i(\nu^{(n-1)})-Q_i(\nu^{(n)}) && \text{for $i\in 2\Z_{\ge 0}$}
\end{aligned}
\end{equation}

\subsection*{Constraints}
\begin{equation}\label{eq:constraint C}
|\nu^{(a)}|=L-\sum_{b=1}^a \la_b \qquad \text{for $1\le a\le n$}
\end{equation}

\subsection*{Algorithm $\delta$}
Set $\ell^{(0)}=0$ and repeat the following process for
$a=1,2,\ldots,n$ or until stopped. Find the minimal index $i\ge
\ell^{(a-1)}$ such that $J^{(a,i)}$ is singular. If no such $i$
exists, set $b=a$ and stop. Otherwise set $\ell^{(a)}=i$ and
continue.

If the process has not stopped continue as follows
for $a=n-1,n-2,\ldots,1$ or until stopped. Set $\lb^{(n)}=\ell^{(n)}$
and reset $\ell^{(n)}=\lb^{(n)}-1$.
If $\ell^{(a)}=\lb^{(a+1)}$ set $\lb^{(a)}=\ell^{(a)}$ and
reset $\ell^{(a)}=\lb^{(a)}-1$. Say case (S) holds.
Otherwise find the minimal index $i\ge \lb^{(a+1)}$ such that $J^{(a,i)}$
is singular. If no such $i$ exists, set $b=\overline{a+1}$.
Otherwise set $\lb^{(a)}=i$ and continue. If the process does not
stop for $a\ge 1$ set $b=\overline{1}$.

Set all undefined $\ell^{(a)}$ and $\lb^{(a)}$ for $1\le a\le n$
to $\infty$.
\subsection*{New RC}
\begin{equation} \label{eq:ch m C1}
m_i^{(a)}(\ntt)=m_i^{(a)}(\nt)+\begin{cases}
 1 & \text{if $i=\ell^{(a)}-1$}\\
 -1 & \text{if $i=\ell^{(a)}$}\\
 1 & \text{if $i=\lb^{(a)}-1$}\\
 -1 & \text{if $i=\lb^{(a)}$}\\
 0 & \text{otherwise} \end{cases}
\end{equation}
If two or more conditions hold then all changes should be performed.

If $a=n$ or case (S) holds for $1\le a<n$ the partition $\tilde{J}^{(a,i)}$ is
obtained from $J^{(a,i)}$ by removing a part of size $P_i^{(a)}(\nt)$ for
$i=\lb^{(a)}$, adding a part of size $P_i^{(a)}(\ntt)$ for $i=\lb^{(a)}-2$,
and leaving it unchanged otherwise.
Otherwise $\tilde{J}^{(a,i)}$ is obtained from $J^{(a,i)}$ by removing a
part of size $P_i^{(a)}(\nt)$ for $i=\ell^{(a)}$ and $i=\lb^{(a)}$,
adding a part of size $P_i^{(a)}(\ntt)$ for $i=\ell^{(a)}-1$ and
$i=\lb^{(a)}-1$, and leaving it unchanged otherwise.

\subsection*{Change in vacancy numbers}
\begin{align}\label{eq:cv C}
P_i^{(a)}(\ntt)=P_i^{(a)}(\nt)-&\chi(\ell^{(a-1)}\le i)+2\chi(\ell^{(a)}\le i)
 -\chi(\ell^{(a+1)}\le i)\\
 -&\chi(\lb^{(a-1)}\le i)+2\chi(\lb^{(a)}\le i)-\chi(\lb^{(a+1)}\le i)\notag\\
\intertext{for $1\le a\le n-1$ and}
P_i^{(n)}(\ntt)=P_i^{(n)}(\nt)
 -&\chi(\ell^{(n-1)}\le i)-\chi(\lb^{(n-1)}\le i)\notag\\
 +&\chi(\ell^{(n)}\le i)+\chi(\lb^{(n)}\le i).\notag
\end{align}

\subsection{Bijection algorithm for type $A_{2n}^{(2)}$}
Recall here that $\gehb=C_n$ and $\gf=B_n$.
\subsection*{Vacancy numbers}
The vacancy numbers are the same as for type $C_n^{(1)}$ \eqref{eq:vac C}
with the only exception that now $i\in\Z_{\ge 0}$ even for $a=n$.

\subsection*{Constraints}
The constraints are the same as for type $C_n^{(1)}$ \eqref{eq:constraint C}.

\subsection*{Algorithm $\delta$}
Set $\ell^{(0)}=0$ and repeat the following process for
$a=1,2,\ldots,n$ or until stopped. Find the minimal index $i\ge
\ell^{(a-1)}$ such that $J^{(a,i)}$ is singular. If no such $i$
exists, set $b=a$ and stop. Otherwise set $\ell^{(a)}=i$ and
continue.

If $\ell^{(n)}=1$ set $b=\phi$ and stop. Otherwise say case (S) holds for
$a=n$ and continue.

If the process has not stopped, set $\lb^{(n)}=\ell^{(n)}$ and reset
$\ell^{(n)}=\lb^{(n)}-1$. Continue as follows for $a=n-1,n-2,\ldots,1$ or
until stopped.
If $\ell^{(a)}=\lb^{(a+1)}$ set $\lb^{(a)}=\ell^{(a)}$ and
reset $\ell^{(a)}=\lb^{(a)}-1$. Say case (S) holds.
Otherwise find the minimal index $i\ge \lb^{(a+1)}$ such that $J^{(a,i)}$
is singular. If no such $i$ exists, set $b=\overline{a+1}$.
Otherwise set $\lb^{(a)}=i$ and continue. If the process does not
stop for $a\ge 1$ set $b=\overline{1}$.

Set all undefined $\ell^{(a)}$ and $\lb^{(a)}$ for $1\le a\le n$
to $\infty$.
\subsection*{New RC}
The configuration changes in the same way as for type $C_n^{(1)}$
\eqref{eq:ch m C1}.

If case (S) holds for $1\le a\le n$ the partition $\tilde{J}^{(a,i)}$ is
obtained from $J^{(a,i)}$ by removing a part of size $P_i^{(a)}(\nt)$ for
$i=\lb^{(a)}$, adding a part of size $P_i^{(a)}(\ntt)$ for $i=\lb^{(a)}-2$,
and leaving it unchanged otherwise.
Otherwise $\tilde{J}^{(a,i)}$ is obtained from $J^{(a,i)}$ by removing a
part of size $P_i^{(a)}(\nt)$ for $i=\ell^{(a)}$ and $i=\lb^{(a)}$,
adding a part of size $P_i^{(a)}(\ntt)$ for $i=\ell^{(a)}-1$ and
$i=\lb^{(a)}-1$, and leaving it unchanged otherwise.

\subsection*{Change in vacancy numbers}
The change in the vacancy numbers is the same as for type $C_n^{(1)}$
\eqref{eq:cv C}.

\subsection{Bijection algorithm for type $A_{2n-1}^{(2)}$}

\subsection*{Vacancy numbers}
\begin{equation}\label{eq:vac A2o}
\begin{split}
P_i^{(a)}(\nt)&=Q_i(\nu^{(a-1)})-2Q_i(\nu^{(a)})+Q_i(\nu^{(a+1)})
 +L\delta_{a,1} \qquad\text{for $1\le a<n-1$}\\
P_i^{(n-1)}(\nt)&=Q_i(\nu^{(n-2)})-2Q_i(\nu^{(n-1)})+2Q_i(\nu^{(n)})\\
P_i^{(n)}(\nt)&=Q_i(\nu^{(n-1)})-2Q_i(\nu^{(n)})
\end{split}
\end{equation}

\subsection*{Constraints}
\begin{equation}\label{eq:constraint A2o}
\begin{aligned}[3]
|\nu^{(a)}|&=L-\sum_{b=1}^a \la_b &&\text{for $1\le a<n$}\\
|\nu^{(n)}|&=\frac{1}{2}(L-\sum_{b=1}^n \la_b)&&
\end{aligned}
\end{equation}

\subsection*{Algorithm $\delta$}
Set $\ell^{(0)}=0$ and repeat the following process for
$a=1,2,\ldots,n$ or until stopped. Find the minimal index $i\ge
\ell^{(a-1)}$ such that $J^{(a,i)}$ is singular. If no such $i$
exists, set $b=a$ and stop. Otherwise set $\ell^{(a)}=i$ and
continue.

If the process has not stopped set $\lb^{(n)}=\ell^{(n)}$ and
continue as follows for $a=n-1,n-2,\ldots,1$ or until stopped.
Find the minimal index $i\ge \lb^{(a+1)}$ such that $J^{(a,i)}$ is singular
(if $i=\ell^{(a)}$ then there need to be two parts of size
$P_i^{(a)}(\nt)$ in $J^{(a,i)}$).
If no such $i$ exists, set $b=\overline{a+1}$ and stop.
If the process did not stop, set $b=\overline{1}$.

Set all yet undefined $\ell^{(a)}$ and $\lb^{(a)}$ to $\infty$.

\subsection*{New RC}
\begin{align} \label{eq:ch m A2o}
m_i^{(a)}(\ntt)=m_i^{(a)}(\nt)+\begin{cases}
 1 & \text{if $i=\ell^{(a)}-1$}\\
 -1 & \text{if $i=\ell^{(a)}$}\\
 1 & \text{if $i=\lb^{(a)}-1$ and $1\le a\le n-1$}\\
 -1 & \text{if $i=\lb^{(a)}$ and $1\le a \le n-1$}\\
 0 & \text{otherwise.} \end{cases}
\end{align}
The partition $\tilde{J}^{(a,i)}$ is obtained from $J^{(a,i)}$ by removing
a part of size $P_i^{(a)}(\nt)$ for $i=\ell^{(a)}$ when $1\le a\le n$
and $i=\lb^{(a)}$ when $1\le a<n$, adding a part of size $P_i^{(a)}(\ntt)$
for $i=\ell^{(a)}-1$ when $1\le a\le n$ and $i=\lb^{(a)}-1$ when
$1\le a<n$, and leaving it unchanged otherwise.

\subsection*{Change in vacancy numbers}
\begin{align}\label{eq:cv A2o}
P_i^{(a)}(\ntt)=P_i^{(a)}(\nt)-&\chi(\ell^{(a-1)}\le i)+2\chi(\ell^{(a)}\le i)
 -\chi(\ell^{(a+1)}\le i)\\
 -&\chi(\lb^{(a-1)}\le i)+2\chi(\lb^{(a)}\le i)-\chi(\lb^{(a+1)}\le i)\notag\\
\intertext{for $1\le a\le n-1$ and}
P_i^{(n)}(\ntt)=P_i^{(n)}(\nt)
 -&\chi(\ell^{(n-1)}\le i)+2\chi(\ell^{(n)}\le i)-\chi(\lb^{(n-1)}\le i).\notag
\end{align}

\subsection{Bijection algorithm for type $D_{n+1}^{(2)}$}
\subsection*{Vacancy numbers}
\begin{equation}\label{eq:vac D2}
\begin{aligned}
P_i^{(a)}(\nt)&=Q_i(\nu^{(a-1)})-2Q_i(\nu^{(a)})+Q_i(\nu^{(a+1)})
 +L\delta_{a,1} && \text{for $1\le a\le n-1$}\\
P_i^{(n)}(\nt)&=2Q_i(\nu^{(n-1)})-2Q_i(\nu^{(n)})&&
\end{aligned}
\end{equation}

\subsection*{Constraints}
The constraints are the same as for type $C_n^{(1)}$ \eqref{eq:constraint C}.

\subsection*{Algorithm $\delta$}
Call a partition quasi-singular if it is not singular and has a
part of size $P_i^{(a)}(\nt)-1$.

Set $\ell^{(0)}=0$ and repeat the following process for
$a=1,2,\ldots,n-1$ or until stopped.
Find the minimal index $i\ge \ell^{(a-1)}$ such that $J^{(a,i)}$ is singular.
If no such $i$ exists, set $b=a$ and stop.
Otherwise set $\ell^{(a)}=i$ and continue.

If the process has not yet stopped, continue as follows. Consider
the following conditions:
\begin{enumerate}
\item[(S)] $J^{(n,i)}$ is singular and $i>1$;
\item[(P)] $J^{(n,i)}$ is singular and $i=1$;
\item[(Q)] $J^{(n,i)}$ is quasi-singular.
\end{enumerate}
Find the minimal index $i\ge \ell^{(n-1)}$ such that one of the
mutually exclusive conditions (S), (P) or (Q) holds. If no such
$i$ exists, set $b=n$ and stop. If (P) holds set
$\ell^{(n)}=i,b=\phi$ and stop. If (S) holds set
$\ell^{(n)}=i-1,\lb^{(n)}=i$, say case (S) holds for $a=n$ and
continue. If (Q) holds set $\ell^{(n)}=i$. Find the minimal $j>i$
such that (S) holds. If no such $j$ exists, set $b=0$ and stop.
Else set $\lb^{(n)}=j$, say case (Q,S) holds and continue.

If the process has not stopped continue in the following fashion
for $a=n-1,n-2,\ldots,1$ or until stopped.
If $\ell^{(a)}=\lb^{(a+1)}$ set $\lb^{(a)}=\ell^{(a)}$ and reset
$\ell^{(a)}=\lb^{(a)}-1$. Say case (S) holds for $a$.
Otherwise find the minimal index $i\ge \lb^{(a+1)}$ such that $J^{(a,i)}$
is singular. If no such $i$ exists, set $b=\overline{a+1}$ and stop.
Otherwise set $\lb^{(a)}=i$ and continue.
If the process did not stop for $a\ge 1$ set $b=\overline{1}$.

Set all undefined $\ell^{(a)}$ and $\lb^{(a)}$ for $1\le a\le n$
to $\infty$.

\subsection*{New RC}
The new configuration $\ntt$ is given by \eqref{eq:ch m C1}.

If case (S) holds for $1\le a\le n$ the partition $\tilde{J}^{(a,i)}$ is
obtained from $J^{(a,i)}$ by removing a part of size $P_i^{(a)}(\nt)$ for
$i=\lb^{(a)}$, adding a part of size $P_i^{(a)}(\ntt)$ for $i=\lb^{(a)}-2$,
and leaving it unchanged otherwise.
If (Q) or (Q,S) holds for $a=n$, then $\tilde{J}^{(n,i)}$ is obtained from
$J^{(n,i)}$ by removing a part of size $P_i^{(n)}(\nt)-1$
(resp. $P_i^{(n)}(\nt)$) for $i=\ell^{(n)}$ (resp. $i=\lb^{(n)}$),
adding a part of size $P_i^{(n)}(\ntt)$ (resp. $P_i^{(n)}(\ntt)-1$)
for $i=\ell^{(n)}-1$ (resp. $i=\lb^{(n)}-1$), and leaving it unchanged
otherwise.
Otherwise $\tilde{J}^{(a,i)}$ is obtained from $J^{(a,i)}$ by removing a
part of size $P_i^{(a)}(\nt)$ for $i=\ell^{(a)}$ and $i=\lb^{(a)}$,
adding a part of size $P_i^{(a)}(\ntt)$ for $i=\ell^{(a)}-1$ and
$i=\lb^{(a)}-1$, and leaving it unchanged otherwise.

\subsection*{Change in vacancy numbers}
\begin{align}\label{eq:cv D2}
P_i^{(a)}(\ntt)=P_i^{(a)}(\nt)-&\chi(\ell^{(a-1)}\le i)+2\chi(\ell^{(a)}\le i)
 -\chi(\ell^{(a+1)}\le i)\\
 -&\chi(\lb^{(a-1)}\le i)+2\chi(\lb^{(a)}\le i)-\chi(\lb^{(a+1)}\le i)\notag\\
\intertext{for $1\le a\le n-1$ and}
P_i^{(n)}(\ntt)=P_i^{(n)}(\nt)
 -&2\chi(\ell^{(n-1)}\le i)+2\chi(\ell^{(n)}\le i)\notag\\
 -&2\chi(\lb^{(n-1)}\le i)+2\chi(\lb^{(n)}\le i).\notag
\end{align}

\subsection{Bijection algorithm for type $A_{2n}^{(2)\dagger}$}
\subsection*{Vacancy numbers}
The vacancy numbers are given by the same formula as for type $C_n^{(1)}$
\eqref{eq:vac C} with the only exception that in this case
$i\in\Z_{\ge 0}$ for all $a\in \I$.

\subsection*{Algorithm $\delta$}
If $a=n$ and $i$ is odd, then $J^{(n,i)}$ is never singular. For
$i$ odd, call $J^{(n,i)}$ quasi-singular if it has a part of size
$P_i^{(n)}(\nt)-1/2$.

Set $\ell^{(0)}=0$ and repeat the following process for
$a=1,2,\ldots,n-1$ or until stopped. Find the minimal index $i\ge
\ell^{(a-1)}$ such that $J^{(a,i)}$ is singular. If no such $i$
exists, set $b=a$ and stop. Otherwise set $\ell^{(a)}=i$ and
continue.

If the process has not yet stopped, continue as follows. Consider
the conditions
\begin{enumerate}
\item[(S)] $i$ is even and $J^{(n,i)}$ is singular;
\item[(Q)] $i$ is odd and $J^{(n,i)}$ is quasi-singular.
\end{enumerate}
Find the minimal index $i\ge \ell^{(n-1)}$ such that one of the
mutually exclusive conditions (S) or (Q) holds. If no such $i$
exists, set $b=n$ and stop. If (S) holds set
$\ell^{(n)}=i-1,\lb^{(n)}=i$, say case (S) holds for $a=n$ and
continue. If (Q) holds set $\ell^{(n)}=i$. Find the minimal $j>i$
such that (S) holds for $j$. If no such $j$ exists, set $b=0$ and
stop. Else set $\lb^{(n)}=j$, say case (Q,S) holds and continue.

If the process has not stopped continue in the following fashion
for $a=n-1,n-2,\ldots,1$ or until stopped. If
$\ell^{(a)}=\lb^{(a+1)}$ set $\lb^{(a)}=\ell^{(a)}$ and reset
$\ell^{(a)}=\lb^{(a)}-1$. Say case (S) holds for $a$. Otherwise
find the minimal index $i\ge \lb^{(a+1)}$ such that $J^{(a,i)}$ is
singular. If no such $i$ exists, set $b=\overline{a+1}$ and stop.
Otherwise set $\lb^{(a)}=i$ and continue. If the process did not
stop for $a\ge 1$ set $b=\overline{1}$.

Set all undefined $\ell^{(a)}$ and $\lb^{(a)}$ for $1\le a\le n$
to $\infty$.

\subsection*{New RC}
The new configuration $\ntt$ is given by \eqref{eq:ch m C1}.

If case (S) holds for $1\le a\le n$ the partition
$\tilde{J}^{(a,i)}$ is obtained from $J^{(a,i)}$ by removing a
part of size $P_i^{(a)}(\nt)$ for $i=\lb^{(a)}$, adding a part of
size $P_i^{(a)}(\ntt)$ for $i=\lb^{(a)}-2$, and leaving it
unchanged otherwise.

If (Q) or (Q,S) holds for $a=n$, then $\tilde{J}^{(n,i)}$ is
obtained from $J^{(n,i)}$ by removing a part of size
$P_i^{(n)}(\nt)-1/2$ for $i=\ell^{(n)}$ (and a part of size
$P_i^{(n)}(\nt)$ for $i=\lb^{(n)}<\infty$), adding a part of size
$P_i^{(n)}(\ntt)$ for $i=\ell^{(n)}-1$ (and a part of size
$P_i^{(n)}(\ntt)-1/2$ for $i=\lb^{(n)}-1<\infty$), and leaving it
unchanged otherwise.

Otherwise $\tilde{J}^{(a,i)}$ is obtained from $J^{(a,i)}$ by
removing a part of size $P_i^{(a)}(\nt)$ for $i=\ell^{(a)}$ and
$i=\lb^{(a)}$, adding a part of size $P_i^{(a)}(\ntt)$ for
$i=\ell^{(a)}-1$ and $i=\lb^{(a)}-1$, and leaving it unchanged
otherwise.

\subsection*{Change in vacancy numbers}
The vacancy numbers $P_i^{(a)}(\nt)$ change as in \eqref{eq:cv C}.

\section{Proof of Theorem \ref{thm:bij}}
\label{sec:proof}

In the following subsections Theorem \ref{thm:bij} is proved
case-by-case for the various root systems. The following notation
is used. Let $(\nt,\Jt)\in\RC(\la,\mu)$, $b=\rk(\nt,\Jt)\in B$,
$\rho=\la-\wt(b)$, and $(\ntt,\Jtt)=\delta(\nt,\Jt)$. There are
three things that must be verified:
\begin{enumerate}
\item[(I)] $\rho$ is dominant and $b$ can be appended to any path in
$\Path(\rho,\mut)$ to give an element of $\Path(\la,\mu)$.
\item[(II)] $(\ntt,\Jtt)\in\RC(\rho,\mut)$ where $B^{\mut}=B^{\otimes
(L-1)}$.
\item[(III)] The conditions of Lemma \ref{lem:stat reduction} are
satisfied.
\end{enumerate}

Parts (I) and (II) show that $\delta$ is well-defined. The proof
that $\delta$ has an inverse, is omitted as it is very similar to
the proof of well-definedness. Part (III) suffices to prove that
$\Phit$ preserves statistics.

For $(\nt,\Jt)\in\RC(\la,\mu)$, define
$\Delta(cc(\nt,\Jt))=cc(\nt,\Jt)-cc(\delta'(\nt,\Jt))$ and
$\Delta^2(cc(\nt,\Jt))=\Delta(cc(\nt,\Jt))-\Delta(cc(\delta'(\nt,\Jt)))$
where $\delta'=\comp\circ\delta\circ\comp$.

\begin{lemma} \label{lem:stat reduction} To prove that
\eqref{eq:cc=D} holds, it suffices to show that it holds for
$L=1$, and that for $L\ge2$ with $\Phit(\nt,\Jt)=b_L\otimes
\dotsm\otimes b_1$, we have
\begin{equation}\label{eq:Deltacc}
  \Delta(cc(\nt,\Jt)) = \frac{t_1^\vee}{a_0^\vee} \alpha^{(1)}_1 -
  \chi(b_L=\phi),
\end{equation}
and
\begin{equation}\label{eq:Delta2cc}
\Hb(b_L\otimes b_{L-1})=
  \frac{t_1^\vee}{a_0^\vee}
  (\alpha^{(1)}_1 -
  \alt^{(1)}_1)-\chi(b_L=\phi)+\chi(b_{L-1}=\phi)
\end{equation}
where $\alpha^{(1)}_1$ and $\alt^{(1)}_1$ are the lengths of the
first columns in $\nu^{(1)}$ and $\tilde{\nu}^{(1)}$ respectively,
and $\delta(\nt,\Jt)=(\ntt,\Jtt)$.
\end{lemma}
\begin{proof} If $L=0$, $\RC(\la,\mu)$ and $\Path(\la,\mu)$ are both
empty unless $\la=0$, in which case $\RC(\la,\mu)$ (resp.
$\Path(\la,\mu)$) is the singleton set containing the empty rigged
configuration (resp. the empty path). Both of these objects have
statistic zero. The case $L=1$ is given by hypothesis. For $L\ge
2$, by the definition \eqref{eq:energy} and \eqref{eq:co} of $\Db$,
\begin{equation} \label{eq:DeltaD}
  \Db(b_L\otimes\dotsm\otimes b_1) -
  \Db(b_{L-1}\otimes\dotsm\otimes b_1) =
  \Hb(b_1 \otimes b^\natural) + \sum_{j=1}^{L-1} \Hb(b_{j+1}\otimes
  b_j).
\end{equation}
Therefore by induction on $L$ it suffices to prove that
$\Delta(cc(\nt,\Jt))$ is given by the right hand side of
\eqref{eq:DeltaD}. By induction and again ``taking the difference"
it suffices to prove that
\begin{equation*}
  \Delta^2(cc(\nt,\Jt)) = \Hb(b_L\otimes b_{L-1}).
\end{equation*}
But this follows from \eqref{eq:Deltacc} and \eqref{eq:Delta2cc}.
\end{proof}

We also need several preliminary lemmas on the convexity and
nonnegativity of the vacancy numbers $P^{(a)}_i(\nt)$.

\begin{lemma}\label{lem:asym}
For large $i$, we have
\begin{equation*}
\begin{split}
P_i^{(a)}(\nt) &= \la_a-\la_{a+1} \qquad\text{for $1\le a<n$}\\
P_i^{(n)}(\nt) &= \begin{cases}
 2\la_n & \text{for $B_n^{(1)}$, $D_{n+1}^{(2)}$}\\
 \la_n & \text{for $C_n^{(1)}$, $A_{2n}^{(2)}$, $A^{(2)\dagger}_{2n}$, $A_{2n-1}^{(2)}$}\\
 \la_{n-1}+\la_n & \text{for $D_n^{(1)}$.} \end{cases}
\end{split}
\end{equation*}
\end{lemma}
\begin{proof}
This follows from the formulas for the vacancy numbers \eqref{eq:vac D},
\eqref{eq:vac B}, \eqref{eq:vac C}, \eqref{eq:vac A2o}, \eqref{eq:vac D2},
the constraints \eqref{eq:constraint D}, \eqref{eq:constraint B},
\eqref{eq:constraint C}, \eqref{eq:constraint A2o},
and the fact that for large $i$, $Q_i(\nu^{(a)})=|\nu^{(a)}|$.
\end{proof}

Direct calculations show that
\begin{align}\label{eq:Pm D}
&\textbf{Type $D_n^{(1)}$}\\
&-P_{i-1}^{(a)}(\nt)+2P_i^{(a)}(\nt)-P_{i+1}^{(a)}(\nt)\notag\\
=&\begin{cases}
m_i^{(a-1)}(\nt)-2m_i^{(a)}(\nt)+m_i^{(a+1)}(\nt)+L\delta_{a,1}\delta_{i,1}
 & \text{for $1\le a\le n-3$}\\
m_i^{(n-3)}(\nt)-2m_i^{(n-2)}(\nt)+m_i^{(n-1)}(\nt)+m_i^{(n)}(\nt)
 & \text{for $a=n-2$}\\
m_i^{(n-2)}(\nt)-2m_i^{(n)}(\nt)
 & \text{for $a=n-1,n$.}
\end{cases}\notag
\end{align}

\begin{align}\label{eq:Pm B}
&\textbf{Type $B_n^{(1)}$}\\
&-P_{i-\f_a}^{(a)}(\nt)+2P_i^{(a)}(\nt)-P_{i+\f_a}^{(a)}(\nt)\notag\\
=&\begin{cases}
m_i^{(a-1)}(\nt)-2m_i^{(a)}(\nt)+m_i^{(a+1)}(\nt)+L\delta_{a,1}\delta_{i,1}
 & \text{for $1\le a\le n-2$}\\
m_i^{(n-2)}(\nt)-2m_i^{(n-1)}(\nt)&\\
\quad
 +2(2m_i^{(n)}(\nt)+m_{i+\frac{1}{2}}^{(n)}(\nt)+m_{i-\frac{1}{2}}^{(n)}(\nt))
 & \text{for $a=n-1$}\\
2m_i^{(n-1)}(\nt)-4m_i^{(n)}(\nt)
 & \text{for $a=n$.}
\end{cases}\notag
\end{align}

\begin{align}\label{eq:Pm C}
&\textbf{Type $C_n^{(1)}$}\\
&-P_{i-\f_a}^{(a)}(\nt)+2P_i^{(a)}(\nt)-P_{i+\f_a}^{(a)}(\nt)\notag\\
=&\begin{cases}
m_i^{(a-1)}(\nt)-2m_i^{(a)}(\nt)+m_i^{(a+1)}(\nt)+L\delta_{a,1}\delta_{i,1}
 & \text{for $1\le a\le n-1$}\\
m_{i-1}^{(n-1)}(\nt)+2m_i^{(n-1)}(\nt)+m_{i+1}^{(n-1)}(\nt)-2m_i^{(n)}(\nt)
 & \text{for $a=n$.}
\end{cases}\notag
\end{align}

\begin{align}\label{eq:Pm A2}
&\textbf{Types $A_{2n}^{(2)}$ and $A^{(2)\dagger}_{2n}$}\\
&-P_{i-1}^{(a)}(\nt)+2P_i^{(a)}(\nt)-P_{i+1}^{(a)}(\nt)\notag\\
=&\begin{cases}
m_i^{(a-1)}(\nt)-2m_i^{(a)}(\nt)+m_i^{(a+1)}(\nt)+L\delta_{a,1}\delta_{i,1}
 & \text{for $1\le a\le n-1$}\\
m_i^{(n-1)}(\nt)-m_i^{(n)}(\nt)
 & \text{for $a=n$.}
\end{cases}\notag
\end{align}

\begin{align}\label{eq:Pm A2o}
&\textbf{Type $A_{2n-1}^{(2)}$}\\
&-P_{i-1}^{(a)}(\nt)+2P_i^{(a)}(\nt)-P_{i+1}^{(a)}(\nt)\notag\\
=&\begin{cases}
m_i^{(a-1)}(\nt)-2m_i^{(a)}(\nt)+m_i^{(a+1)}(\nt)+L\delta_{a,1}\delta_{i,1}
 & \text{for $1\le a<n-1$}\\
m_i^{(n-2)}(\nt)-2m_i^{(n-1)}(\nt)+2m_i^{(n)}(\nt)
 & \text{for $a=n-1$}\\
m_i^{(n-1)}(\nt)-2m_i^{(n)}(\nt)
 & \text{for $a=n$.}
\end{cases}\notag
\end{align}

\begin{align}\label{eq:Pm D2}
&\textbf{Type $D_{n+1}^{(2)}$}\\
&-P_{i-1}^{(a)}(\nt)+2P_i^{(a)}(\nt)-P_{i+1}^{(a)}(\nt)\notag\\
=&\begin{cases}
m_i^{(a-1)}(\nt)-2m_i^{(a)}(\nt)+m_i^{(a+1)}(\nt)+L\delta_{a,1}\delta_{i,1}
 & \text{for $1\le a\le n-1$}\\
2m_i^{(n-1)}(\nt)-2m_i^{(n)}(\nt)
 & \text{for $a=n$.}
\end{cases}\notag
\end{align}

In particular these equations imply the convexity condition
\begin{equation}\label{eq:convex}
P_i^{(a)}(\nt)\ge \frac{1}{2}(P_{i-\f_a}^{(a)}(\nt)+P_{i+\f_a}^{(a)}(\nt))
\qquad \text{if $m_i^{(a)}(\nt)=0$.}
\end{equation}

\begin{lemma}\label{lem:equiv}
Let $\nt$ be a configuration in $C(\la,\mu)$. The following
are equivalent:
\begin{enumerate}
\item $P_i^{(a)}(\nt)\ge 0$ for all $i\in\f_a\Z_{>0}$, $a\in\I$;
\item $P_i^{(a)}(\nt)\ge 0$ for all $i\in\f_a\Z_{>0}$, $a\in\I$ such that
$m_i^{(a)}(\nt)>0$.
\end{enumerate}
\end{lemma}
\begin{proof}
This follows immediately from Lemma \ref{lem:asym} and the convexity
condition \eqref{eq:convex}. (See also \cite[Lemma 10]{KS}).
\end{proof}

\subsection{Proof for type $D_n^{(1)}$}

\begin{proof}[Proof of (I) for $D_n^{(1)}$]
Here it suffices to show that $\rho$
satisfies \eqref{eq:dominant D}. Suppose not. If $b=k$ with $1\le
k\le n$ then
\begin{equation*}
\begin{split}
&\text{(a) $\la_k=\la_{k+1}$ if $1\le k\le n-2$}\\
&\text{(b) $\la_{n-1}=|\la_n|$ if $k=n-1$}\\
&\text{(c) $\la_{n-1}=-\la_n$ if $k=n$.}
\end{split}
\end{equation*}
In case (a) we have $P_i^{(k)}(\nt)=0$ for large $i$ by Lemma
\ref{lem:asym}. Let $\ell$ be the largest part in $\nu^{(k)}$. By
convexity this implies $P_i^{(k)}(\nt)=0$ for all $i\ge \ell$.
Equation \eqref{eq:Pm D} in turn yields $m_i^{(k-1)}(\nt)=0$ for
all $i>\ell$ so that $1\le \ell^{(k-1)}\le \ell$. But this is a
contradiction since there is a singular string of length $\ell$ in
$(\nt,\Jt)^{(k)}$ since $P_\ell^{(k)}(\nt)=0$ and
$m_\ell^{(k)}(\nt)>0$ so that we would have $\rk(\nt,\Jt)>k$. In
case (b) let us first assume that $\la_{n-1}=\la_n$. Then for
large $i$, $P_i^{(n-1)}(\nt)=0$ and by convexity
$P_i^{(n-1)}(\nt)=0$ for $i\ge \ell$ where $\ell$ is the largest
part in $\nu^{(n-1)}$. By \eqref{eq:Pm D} we have
$m_i^{(n-2)}(\nt)=0$ for $i>\ell$. Hence $1\le \ell^{(n-2)}\le
\ell$ which yields a contradiction since there is a singular
string of length $\ell$ in $(\nt,\Jt)^{(n-1)}$ so that
$\rk(\nt,\Jt)\neq n-1$. If $\la_{n-1}=-\la_n$ the same argument
goes through with $n-1$ replaced by $n$. The case (c) is analogous
to the second part of case (b).

Now suppose $b=\overline{k}$ for some $1\le k\le n$. We show again
that $\rho$ not dominant will yield a contradiction. If $\rho$ is
not dominant one of the following has to be true:
\begin{equation*}
\begin{split}
&\text{(d) $\la_k=\la_{k-1}$ if $2\le k\le n-1$}\\
&\text{(e) $\la_n=\la_{n-1}$ if $k=n$.}
\end{split}
\end{equation*}
Case (e) is analogous to case (b). In case (d) some caution is
in order. By lemma \ref{lem:asym} and convexity \eqref{eq:convex}
we have $P_i^{(k-1)}(\nt)=0$ for $i\ge \ell$ where $\ell$ is the largest
part in $\nu^{(k-1)}$. By \eqref{eq:Pm D} it follows that
$m_i^{(k)}(\nt)=0$ for $i>\ell$. Hence $\lb^{(k)}\le \ell$.
Since $P_\ell^{(k-1)}(\nt)=0$ and $m_\ell^{(k-1)}(\nt)>0$ there is
a singular string of length $\ell$ in $(\nt,\Jt)^{(k-1)}$.
Hence $\lb^{(k-1)}\le \ell$ unless $\ell^{(k-1)}=\ell$ and
$m_\ell^{(k-1)}(\nt)=1$. We will show that the latter case cannot
occur. Equation \eqref{eq:Pm D} with $a=k-1$ and $i=\ell$ implies that
$P_{\ell-1}^{(k-1)}=0$ and $m_\ell^{(k-2)}(\nt)=0$ since
by assumption $\ell^{(k-1)}=\ell^{(k)}=\lb^{(k)}=\ell$ and hence
$m_\ell^{(k)}(\nt)\ge 2$ (or $m_\ell^{(n-1)}(\nt)\ge 1$ and
$m_\ell^{(n)}(\nt)\ge 1$ for $k=n-1$). However this implies that
$m_{\ell-1}^{(n-2)}(\nt)=0$ since otherwise $\ell^{(k-1)}\le \ell-1$
and not $\ell$ since there is a singular string of length $\ell-1$
in $(\nt,\Jt)^{(k-1)}$. Now by induction on $i=\ell-1,\ell-2,\ldots,1$
it follows from \eqref{eq:Pm D} at $a=k-1$ that $P_i^{(k-1)}(\nt)=
m_i^{(k-2)}(\nt)=m_i^{(k-1)}(\nt)=0$. However, this means in particular
that $m_i^{(k-2)}(\nt)=0$ for all $1\le i\le \ell$ so that
$\ell^{(k-2)}>\ell$ which contradicts $\ell^{(k-1)}=\ell$.
\end{proof}

\begin{proof}[Proof of (II) for $D_n^{(1)}$]
Denote by $\Jm^{(a,i)}(\nt,\Jt)$ the biggest part in $J^{(a,i)}$.
To prove admissibility of $(\ntt,\Jtt)$ we need to show for all
$i\ge 1, 1\le a\le n$ that
\begin{equation}\label{eq:ineq D}
0\le \Jm^{(a,i)}(\ntt,\Jtt)\le P_i^{(a)}(\ntt).
\end{equation}
Fix $a\ge 1$. Only one string of size $\ell^{(a)}$ and one string of size
$\lb^{(a)}$ change in the transformation $(\nt,\Jt)^{(a)}\to (\ntt,\Jtt)^{(a)}$.
Hence
\begin{align*}
&\Jm^{(a,i)}(\ntt,\Jtt)=P_i^{(a)}(\ntt)
 && \text{for $i=\ell^{(a)}-1$ and $i=\lb^{(a)}-1$}\\
&0 \le \Jm^{(a,i)}(\ntt,\Jtt)\le \Jm^{(a,i)}(\nt,\Jt)
 && \text{else.}
\end{align*}
Hence by \eqref{eq:cv D} the inequality \eqref{eq:ineq D} can only be violated
when $\ell^{(a-1)}\le i<\ell^{(a)}$ or $\lb^{(a+1)}\le i<\lb^{(a)}$
where $\lb^{(n-1)}=\max(\ell^{(n-1)},\ell^{(n)})$. By the construction of
$\ell^{(a)}$ and $\lb^{(a)}$ there are no singular strings of length
$i$ in $(\nt,\Jt)^{(a)}$ for $\ell^{(a-1)}\le i<\ell^{(a)}$ or
$\lb^{(a+1)}\le i<\lb^{(a)}$. This means that
$\Jm^{(a,i)}(\nt,\Jt)\le P_i^{(a)}(\nt)-1$ if $i$ occurs as a part in
$\nu^{(a)}$, that is $m_i^{(a)}(\nt)>0$. Hence \eqref{eq:ineq D} is fulfilled
for these $i$. It remains to prove that $P_i^{(a)}(\ntt)\ge 0$
for all $i$ such that $m_i^{(a)}(\nt)=0$ and $\ell^{(a-1)}\le i<\ell^{(a)}$
or $\lb^{(a+1)}\le i<\lb^{(a)}$. Note that $m_i^{(a)}(\ntt)=0$
if $m_i^{(a)}(\nt)=0$ for $\ell^{(a-1)}\le i<\ell^{(a)}-1$ or
$\lb^{(a+1)}\le i<\lb^{(a)}-1$. Hence by lemma \ref{lem:equiv}
it suffices to prove \eqref{eq:ineq D} for all $a$ and $i$ such that
$m_i^{(a)}(\ntt)>0$. Therefore the only remaining case for which
\eqref{eq:ineq D} might be violated occurs when
\begin{equation*}
\begin{split}
&\text{For $1\le a\le n-2$:}\\
&\text{$m_{\ell-1}^{(a)}(\nt)=0$, $P_{\ell-1}^{(a)}(\nt)=0$,
$\ell^{(a-1)}<\ell$ (resp. $\lb^{(a+1)}<\ell$)}\\
&\qquad \text{and $\ell$ finite where $\ell=\ell^{(a)}$ (resp.
 $\ell=\lb^{(a)}$)}\\
&\text{For $a=n-1,n$:}\\
&\text{$m_{\ell-1}^{(a)}(\nt)=0$, $P_{\ell-1}^{(a)}(\nt)=0$,
$\ell^{(n-2)}<\ell$}\\
&\qquad \text{and $\ell$ finite where $\ell=\ell^{(a)}$.}
\end{split}
\end{equation*}
We show that these conditions cannot be met simultaneously.
Let $p<\ell$ be maximal such that $m_p^{(a)}(\nt)>0$; if no such
$p$ exists set $p=0$. By \eqref{eq:convex} $P_{\ell-1}^{(a)}(\nt)=0$
is only possible if $P_i^{(a)}(\nt)=0$ for all $p\le i\le \ell$.
By \eqref{eq:Pm D} we find that $m_i^{(a-1)}(\nt)=0$ (resp.
$m_i^{(a+1)}(\nt)=0$) for $p<i<\ell$. Since $\ell^{(a-1)}<\ell$
(resp. $\lb^{(a+1)}<\ell$) this implies that $\ell^{(a-1)}\le p$
(resp. $\lb^{(a+1)}\le p$).
If $p=0$ this contradicts the condition $\ell^{(a-1)}\ge 1$ (resp.
$\lb^{(a+1)}\ge 1$). Hence assume that $p>0$. Since $P_p^{(a)}(\nt)=0$
and $m_p^{(a)}(\nt)>0$ there is a singular string of length $p$ in
$(\nt,\Jt)^{(a)}$ and therefore $\ell^{(a)}=p$ (resp. $\lb^{(a)}=p$).
However, this contradicts $p<\ell$. This concludes the proof that
$(\ntt,\Jtt)$ is well-defined.
\end{proof}

\begin{proof}[Proof of (III) for $D_n^{(1)}$] Here $b^\natural=1$,
$\Hb(b\otimes b')=0$ if $b\le b'$, $\Hb(b\otimes b')=1$ if
$b\otimes b'= n\otimes \overline{n}, \overline{n}\otimes n$ or
$b>b'$ where $b\neq \overline{1}$, $b'\neq 1$, and
$H(\overline{1}\otimes 1)=2$.

If $L=1$ then the path is $1$, the rigged configuration is empty,
and both sides of \eqref{eq:cc=D} are zero.

Here \eqref{eq:Deltacc} and \eqref{eq:Delta2cc} are given by
\begin{align}\label{eq:change cc D1}
\Delta(cc(\nt,\Jt))&=\alpha_1^{(1)} \\
\label{eq:dc D D1} \Hb(b_L\otimes b_{L-1})
 &=\chi(\ell^{(1)}=1)+\chi(\lb^{(1)}=1)
\end{align}
where $\ell^{(i)}$ and $\lb^{(i)}$ are determined by the algorithm
$\delta$.

Let $\lt^{(a)}$ and $\tilde{\overline{\ell}}^{(a)}$ be the length
of the selected strings defined by the algorithm $\delta$ on
$(\ntt,\Jtt)= \delta(\nt,\Jt)$. To check \eqref{eq:dc D D1} note
that if $\ell^{(1)}=1$ it follows from \eqref{eq:cv D} that
$\lt^{(a)}\ge \ell^{(a+1)}$ for $1\le a\le n-2$. Hence if $b_L\le n-1$
then $b_{L-1}<b_L$ and both sides of \eqref{eq:dc D D1} yield 1. If $b_L=n$
then $b_{L-1}\le n-1$ or $b_{L-1}=\overline{n}$ and both sides of
\eqref{eq:dc D D1} are 1. Similarly, if $b_L=\overline{n}$ then
$b_{L-1}\le n$ and both sides if \eqref{eq:dc D D1} are 1.
Finally, if $b_L\ge \overline{n-1}$ then
$\tilde{\overline{\ell}}^{(a)}\ge \lb^{(a-1)}$ and $b_{L-1}<b_L$. If
$b_L<\overline{1}$ then both sides of \eqref{eq:dc D D1} are 1. If
$b_L=\overline{1}$ and $\lb^{(1)}=1$ then $\lt^{(1)}=\infty$ and
hence $b_{L-1}=1$. In this case both sides of \eqref{eq:dc D D1}
are 2. If $b_L=\overline{1}$ and $\lb^{(1)}>1$ then there is a
singular string in $(\ntt,\Jtt)^{(1)}$ so that $b_{L-1}>1$. In
this case both sides of \eqref{eq:dc D D1} are 1. If
$\ell^{(1)}>1$ then $\lt^{(a)}<\ell^{(a)}$ for $1\le a\le n-2$ and
the cases can be checked in a similar fashion as before.

To prove \eqref{eq:change cc D1}, by \eqref{eq:cc} and \eqref{eq:ch m D1}
we have
\begin{equation*}
\begin{split}
&cc(\ntt)=\frac{1}{2}
 \sum_{i,j\ge 1} \sum_{a,b=1}^n \min(i,j) (\alpha_a|\alpha_b)\\
 &\times \left(m_i^{(a)}-\delta_{i,\ell^{(a)}}+\delta_{i,\ell^{(a)}-1}
 -\chi(a\le n-2)(\delta_{i,\lb^{(a)}}-\delta_{i,\lb^{(a)}-1})\right)\\
 &\times\left(m_j^{(b)}-\delta_{j,\ell^{(b)}}+\delta_{j,\ell^{(b)}-1}
 -\chi(b\le n-2)(\delta_{j,\lb^{(b)}}-\delta_{j,\lb^{(b)}-1})\right).
\end{split}
\end{equation*}
Applying the data for $D_n$ and using \eqref{eq:cv D}, a
tedious but straightforward calculation yields
\begin{multline*}
\Delta cc(\nt)=\sum_{a=1}^n \sum_{i\ge 1}
 \left(P_i^{(a)}(\ntt)-P_i^{(a)}(\nt)\right)\\
 \times\left( m_i^{(a)}-\delta_{i,\ell^{(a)}}-\chi(a\le n-2)\delta_{i,\lb^{(a)}})\right)
 +\sum_{i\ge 1}m_i^{(1)}.
\end{multline*}
For $\Delta |\Jt|$ we obtain from the algorithm $\delta'$
\begin{equation*}
\Delta |\Jt|=\sum_{a=1}^n \sum_{i\ge 1}
 \left(P_i^{(a)}(\nt)-P_i^{(a)}(\ntt)\right)
 \left(m_i^{(a)}-\delta_{i,\ell^{(a)}}-\chi(a\le n-2)\delta_{i,\lb^{(a)}})\right).
\end{equation*}
Hence altogether, using $\sum_{i\ge 1}m_i^{(1)}=\alpha_1^{(1)}$,
we obtain \eqref{eq:change cc D1}.
\end{proof}

\subsection{Proof for type $B_n^{(1)}$}

\begin{proof}[Proof of (I) for $B_n^{(1)}$]
Let us assume that either $\rho$ is not dominant, or that $b=0$
(so that $\rho=\la$) and $\la_n=0$. For $b=k$ with $1\le k<n$ the
proof that this cannot happen is the same as for type $D_n^{(1)}$.
Now assume that $b=n$ and $\la_n=0$. Then $P_i^{(n)}(\nt)=0$ for
$i\ge \ell$ where $\ell$ is the largest part in $\nu^{(n)}$ by
Lemma \ref{lem:asym} and \eqref{eq:convex}. By \eqref{eq:Pm B}
with $a=n$ we find that $m_i^{(n-1)}(\nt)=0$ for $i>\ell$, so that
$\ell^{(n-1)}\le \ell$. But there is a singular string of length
$\ell$ in $(\nt,\Jt)^{(n)}$ which contradicts $\ell^{(n)}=\infty$.
Next assume that $b=0$ and $\la_n=0$. By the same arguments as in
the previous case $\ell^{(n-1)}\le \ell$. But there is a singular
string of length $\ell$ in $(\nt,\Jt)^{(n)}$ since
$m_\ell^{(n)}(\nt)>0$ and $P_\ell^{(n)}(\nt)=0$.\ Since (Q) must
hold for $b=0$, there must be a singular string at
$\ell^{(n-1)}-\frac{1}{2}$ or a quasisingular string at
$\ell^{(n-1)}\le i<\ell$. But then (S) holds for $\ell$ which
contradicts $b=0$. The case $b=\overline{k}$ with $1\le k\le n$ is
the same as for type $D_n^{(1)}$.
\end{proof}

\begin{proof}[Proof of (II) for $B_n^{(1)}$]
Denote by $\Jm^{(a,i)}(\nt,\Jt)$ the biggest part in $J^{(a,i)}$.
To prove admissibility of $(\ntt,\Jtt)$ we need to show for all
$i\ge 1, 1\le a\le n$ that
\begin{equation}\label{eq:ineq B}
0\le \Jm^{(a,i)}(\ntt,\Jtt)\le P_i^{(a)}(\ntt).
\end{equation}
Up to small alterations, the proof of \eqref{eq:ineq B} for $1\le a<n$
is the same as for type $D_n^{(1)}$. Let us assume that $a=n$.
Only one string of size $\ell^{(n)}$ and one string of size
$\lb^{(n)}$ change in the transformation $(\nt,\Jt)^{(n)}\to
(\ntt,\Jtt)^{(n)}$. Hence for the different cases:
\begin{equation*}
\begin{aligned}
(S) && \Jm^{(n,i)}(\ntt,\Jtt)&=P_i^{(n)}(\ntt) && \text{for $i=\lb^{(n)}-1$}\\
    &&0 \le \Jm^{(n,i)}(\ntt,\Jtt)&\le \Jm^{(n,i)}(\nt,\Jt) &&
    \text{else} \\
(Q) && \Jm^{(n,i)}(\ntt,\Jtt)&=P_i^{(n)}(\ntt)
 && \text{for $i=\ell^{(n)}-1/2$}\\
    &&0 \le \Jm^{(n,i)}(\ntt,\Jtt)&\le \Jm^{(n,i)}(\nt,\Jt)
 && \text{else}\\
(Q,S) && \Jm^{(n,i)}(\ntt,\Jtt)&=
   P_i^{(n)}(\ntt) && \text{for $i=\ell^{(n)}-1/2$}\\
   && \Jm^{(n,i)}(\ntt,\Jtt)&= P_i^{(n)}(\ntt) &&
   \text{for $i=\lb^{(n)}-1/2$, $\lb^{(n)}=\lb^{(n-1)}$}\\
   && \Jm^{(n,i)}(\ntt,\Jtt)&=P_i^{(n)}(\ntt)-1 && \text{for $i=\lb^{(n)}-1/2$,
     $\lb^{(n)}<\lb^{(n-1)}$}
 \\
    &&0 \le \Jm^{(n,i)}(\ntt,\Jtt)&\le \Jm^{(n,i)}(\nt,\Jt)
 && \text{else}
\end{aligned}
\end{equation*}

Let us first assume that (S) holds:\newline
By the definition of $\ell^{(n)}$ and $\lb^{(n)}$ there is no singular
string at $\ell^{(n-1)}-\frac{1}{2}$ and no singular or quasisingular
string of length $\ell^{(n-1)}\le i<\lb^{(n)}=\ell$.
Hence, if $m_i^{(n)}(\nt)>0$, we have $\Jm^{(n,i)}(\nt,\Jt)
\le P_i^{(n)}(\nt)-2$ for $\ell^{(n-1)}\le i<\ell$ and
$\Jm^{(n,i)}(\nt,\Jt)\le P_i^{(n)}(\nt)-1$ for $i=\ell^{(n-1)}-\frac{1}{2}$.
Hence \eqref{eq:ineq B} holds if $m_i^{(n)}(\nt)>0$.
By lemma \ref{lem:equiv}, \eqref{eq:ineq B} can only be violated if
\begin{equation*}
m_{\ell-1}^{(n)}(\nt)=0,\quad P_{\ell-1}^{(n)}(\nt)=\text{$0$ or $1$},\quad
\ell^{(n-1)}\le \ell-1,\quad \text{$\ell$ finite.}
\end{equation*}
The case $P_{\ell-1}^{(n)}(\nt)=0$ is the same as before.
Hence assume that $P_{\ell-1}^{(n)}(\nt)=1$. If
$m_{\ell-\frac{1}{2}}(\nt)=0$, then by \eqref{eq:Pm B} and \eqref{eq:convex}
$P_i^{(n)}(\nt)=1$ for $p\le i<\ell$ where $p<\ell$ is maximal such
that $m_p^{(n)}(\nt)>0$. By \eqref{eq:Pm B} we also have
$m_i^{(n-1)}(\nt)=0$ for $p<i<\ell$ so that $\ell^{(n-1)}\le p$.
But since $P_p^{(n)}(\nt)=1$ and $m_p^{(n)}(\nt)>0$ there is a
(quasi)singular string of length $p$ in $(\nt,\Jt)^{(n)}$ which
contradicts $p<\ell$.
If $m_{\ell-\frac{1}{2}}^{(n)}(\nt)>0$, then
$P_{\ell-\frac{1}{2}}^{(n)}(\nt)\ge 2$ since otherwise there would be
a (quasi)singular string of length $\ell-\frac{1}{2}$ in $(\nt,\Jt)^{(n)}$.
By convexity \eqref{eq:convex} and \eqref{eq:Pm B} this implies
$P_{\ell-\frac{3}{2}}^{(n)}(\nt)=0$ and $m_{\ell-\frac{3}{2}}(\nt)>0$.
Since $\ell^{(n-1)}\le \ell-1$, (Q) would hold for $\ell-\frac{3}{2}$
which contradicts our assumptions.

One more problem might occur when $\lb^{(n-1)}=\lb^{(n)}=\ell$ and
$m_{\ell-\frac{1}{2}}(\nt)>0$,
$\Jm^{(n,\ell-\frac{1}{2})}(\nt,\Jt)
=P_{\ell-\frac{1}{2}}^{(n)}(\nt)$. But in this case there is a
singular string of length $\ell-\frac{1}{2}$ in $(\nt,\Jt)^{(n)}$
which contradicts $\lb^{(n)}=\ell$.

Now assume that (Q) holds:\newline
In this case $\lb^{(n-1)}=\lb^{(n)}=\infty$. By similar arguments
as before \eqref{eq:ineq B} can only be violated if
\begin{align*}
&m_{\ell-\frac{1}{2}}^{(n)}(\nt)=0,\quad
&P_{\ell-\frac{1}{2}}^{(n)}(\nt)=0,\quad
&\ell^{(n-1)}-\frac{1}{2}<\ell^{(n)}<\ell^{(n)}=\ell,\quad
&\text{$\ell$ finite}\\
\text{or}\quad &m_{\ell-\frac{1}{2}}^{(n)}(\nt)=0,\quad
&P_{\ell-\frac{1}{2}}^{(n)}(\nt)=1,\quad
&\ell^{(n-1)}<\ell^{(n)}<\ell^{(n)}=\ell,\quad &\text{$\ell$
finite.}
\end{align*}
Since (Q) holds, we must have $P_\ell^{(n)}(\nt)\ge 1$. Hence by
convexity \eqref{eq:convex} it follows that
$P_i^{(n)}(\nt)=1$ for $p\le i<\ell$ where $p<\ell$ is maximal such that
$m_p^{(n)}(\nt)>0$.
Equation \eqref{eq:Pm B} implies that $m_i^{(n-1)}(\nt)=0$ for $p<i<\ell$
so that $\ell^{(n-1)}\le p$. But there is a (quasi)singular string
of length $p$ in $(\nt,\Jt)^{(n)}$ which contradicts $p<\ell$.

Finally assume that (Q,S) holds:\newline
For $i<\ell^{(n)}$ the same arguments hold as for case (Q).
Since by definition there are no singular strings of length
$\ell^{(n-1)}<i<\lb^{(n)}=\ell$ in $(\nt,\Jt)^{(n)}$, case (Q) holds
for $i=\ell^{(n-1)}$ and $m_{\ell^{(n)}}^{(n)}(\nt)>0$, the only problem
occurs when
\begin{equation*}
m_{\ell-\frac{1}{2}}^{(n)}(\nt)=0,\quad P_{\ell-\frac{1}{2}}^{(n)}(\nt)=0,
\quad \ell^{(n)}+1\le \ell,\quad \text{$\ell$ finite}.
\end{equation*}
If $p<\ell$ is maximal such that $m_p^{(n)}(\nt)>0$, then by \eqref{eq:Pm B}
and \eqref{eq:convex} $P_i^{(n)}(\nt)=0$ for $p\le i\le \ell$.
Since $m_{\ell^{(n)}}^{(n)}(\nt)>0$ and $P_{\ell^{(n)}}^{(n)}(\nt)>0$
we must have $\ell^{(n)}<p$. But then there is a singular string
of length $p$ in $(\nt,\Jt)^{(n)}$ which contradicts $\lb^{(n)}=\ell$.
\end{proof}

\begin{proof}[Proof of (III) for $B_n^{(1)}$] Here
$b^\natural=\overline{1}$. Note that $\Hb(b\otimes b')=0$ if $b\le
b'$ and $b\otimes b'\neq 0\otimes 0$, $\Hb(b\otimes b')=2$ if
$b\otimes b'=\overline{1}\otimes 1$, and $H(b\otimes b')=1$
otherwise.

If $L=1$ then the path is $1$, the rigged configuration is empty,
and both sides of \eqref{eq:cc=D} are zero.

Here \eqref{eq:Deltacc} and \eqref{eq:Delta2cc} are given by
\begin{align}\label{eq:change cc B1}
\Delta(cc(\nt,\Jt))&=\alpha_1^{(1)} \\
\label{eq:dc D B1} \Hb(b_L\otimes b_{L-1})
 &=\chi(\ell^{(1)}=1)+\chi(\lb^{(1)}=1).
\end{align}
where $\ell^{(i)}$ and $\lb^{(i)}$ are determined by the algorithm
$\delta$.

Let $\lt^{(a)}$ and $\tilde{\overline{\ell}}^{(a)}$ be the length
of the selected strings defined by the algorithm $\delta$ on
$(\ntt,\Jtt)= \delta(\nt,\Jt)$. To check \eqref{eq:dc D B1} note
that if $\ell^{(1)}=1$ it follows that $\lt^{(a)}\ge \ell^{(a+1)}$
for $1\le a\le n-1$. Hence if $b_L\le n$ then $b_{L-1}<b_L$ and
both sides of \eqref{eq:dc D B1} yield 1. If $b_L=0$ then
$b_{L-1}\le b_L$ by \eqref{eq:cv B} and both sides of \eqref{eq:dc
D D1} are 1. If $\overline{n}\le b_L<\overline{1}$, then
$b_{L-1}<b_L$ by \eqref{eq:cv B} and both sides of \eqref{eq:dc D
B1} are 1. If $b_L=\overline{1}$ and $\lb^{(1)}=1$, then
$b_{L-1}=1$ by \eqref{eq:cv B}. Hence both sides of \eqref{eq:dc D
B1} yield 2. Finally, if $b_L=\overline{1}$ and $\lb^{(1)}>1$,
then there exists a singular string of length $\lb^{(1)}-1$ in
$(\ntt,\Jtt)$ so that $b_{L-1}\neq 1$. Hence both sides of
\eqref{eq:dc D B1} are 1. If $\ell^{(1)}>1$ then
$\lt^{(a)}<\ell^{(a)}$ for $1\le a\le n-2$ and the cases can be
checked in a similar fashion as before.

To prove \eqref{eq:change cc B1}, by \eqref{eq:cc} and
\eqref{eq:ch m B1} we have
\begin{equation*}
\begin{split}
&cc(\ntt)=\frac{1}{2}
 \sum_{i,j\ge 1} \sum_{a,b=1}^n \min(t_b i,t_a j)(\alpha_a|\alpha_b)\\
 &\times \left(m_i^{(a)}-\delta_{i,\frac{\ell^{(a)}}{\f_a}}
  +\delta_{i,\frac{\ell^{(a)}}{\f_a}-1}
  -\delta_{i,\frac{\lb^{(a)}}{\f_a}}
  +\delta_{i,\frac{\lb^{(a)}}{\f_a}-1}\right)\\
 &\times\left(m_j^{(b)}-\delta_{j,\frac{\ell^{(b)}}{\f_b}}
  +\delta_{j,\frac{\ell^{(b)}}{\f_b}-1}
  -\delta_{j,\frac{\lb^{(b)}}{\f_b}}
  +\delta_{j,\frac{\lb^{(b)}}{\f_b}-1}\right).
\end{split}
\end{equation*}
Applying the data for $B_n$ and using \eqref{eq:cv B}, a tedious
but straightforward calculation yields
\begin{multline*}
\Delta cc(\nt)=\sum_{a=1}^n \sum_{i\ge 1}
 \left(P_{\f_a i}^{(a)}(\ntt)-P_{\f_a i}^{(a)}(\nt)\right)
 \left(m_i^{(a)}-\delta_{\f_a i,\ell^{(a)}}-\delta_{\f_a i,\lb^{(a)}}\right)
 +\sum_{i\ge 1}m_i^{(1)}\\
 -\chi(\ell^{(n)}=\ell^{(n-1)}-\frac{1}{2})+\chi(\lb^{(n)}=\lb^{(n-1)})
 +\chi(\lb^{(n)}=\infty)\chi(\ell^{(n)}<\infty).
\end{multline*}
For $\Delta |\Jt|$ we obtain from the algorithm $\delta'$
\begin{multline*}
\Delta |\Jt|=\sum_{a=1}^n \sum_{i\ge 1}
 \left(P_{\f_a i}^{(a)}(\nt)-P_{\f_a i}^{(a)}(\ntt)\right)
 \left(m_i^{(a)}-\delta_{\f_a i,\ell^{(a)}}-\delta_{\f_a i,\lb^{(a)}}\right)\\
 +\chi(\ell^{(n)}=\ell^{(n-1)}-\frac{1}{2})-\chi(\lb^{(n)}=\lb^{(n-1)})
 -\chi(\lb^{(n)}=\infty)\chi(\ell^{(n)}<\infty),
\end{multline*}
where the last three terms come from the fact that for $n$-th rigged
partition singular strings can be transformed into quasisingular strings
and vice versa.
Hence altogether, using $\sum_{i\ge 1}m_i^{(1)}=\alpha_1^{(1)}$,
we obtain \eqref{eq:change cc B1}.
\end{proof}

\subsection{Proof for type $C_n^{(1)}$}

\begin{proof}[Proof of (I) for $C_n^{(1)}$]
If $b=k$ with $1\le k<n$ the proof that
$\rho$ is dominant is analogous to type $D_n^{(1)}$. For $b=n$ a
problem occurs if $\la_n=0$. In this case $P_i^{(n)}(\nt)=0$ for
$i\ge \ell$ where $\ell$ is the largest part in $\nu^{(n)}$ by
Lemma \ref{lem:asym} and \eqref{eq:convex}. By \eqref{eq:Pm C}
this implies $m_i^{(n-1)}(\nt)=0$ for $i>\ell$. Hence
$\ell^{(n-1)}\le \ell$. But there is a singular string of length
$\ell$ in $(\nt,\Jt)^{(n)}$ which contradicts $\ell^{(n)}=\infty$.
If $k=\overline{n}$ a problem occurs if $\la_n=\la_{n-1}$. In this
case $P_i^{(n-1)}(\nt)=0$ for $i\ge \ell$ where $\ell$ is the
largest part in $\nu^{(n-1)}$. By \eqref{eq:Pm C},
$m_i^{(n)}(\nt)=0$ for $i>\ell$. Hence $\lb^{(n)}\le \ell$. But
there is a singular string of size $\ell$ in $(\nt,\Jt)^{(n-1)}$
(this also works if $\ell^{(n-1)}=\lb^{(n)}=\ell$) which
contradicts $\lb^{(n-1)}=\infty$.
\end{proof}

\begin{proof}[Proof of (II) for $C_n^{(1)}$]
We show that $(\ntt,\Jtt)\in\RC(\rho,\mut)$. We use the same
notation and set-up as in type $D_n^{(1)}$. Then
\begin{align*}
&\Jm^{(a,i)}(\ntt,\Jtt)=P_i^{(a)}(\ntt)
 && \text{for $i=\ell^{(a)}-1$ and $i=\lb^{(a)}-1$}\\
& && \text{or $i=\lb^{(a)}-2$ if $\ell^{(a)}=\lb^{(a+1)}$}\\
&0 \le \Jm^{(a,i)}(\ntt,\Jtt)\le \Jm^{(a,i)}(\nt,\Jt)
 && \text{else.}
\end{align*}
The proof that $0\le \Jm^{(a,i)}(\ntt,\Jtt)\le P_i^{(a)}(\ntt)$
for $1\le a<n$ is the same as usual if $\ell^{(a)}\neq \lb^{(a+1)}$.
If $\ell^{(a)}=\lb^{(a+1)}$, by \eqref{eq:cv C} the only problem
occurs if
\begin{equation*}
m_{\ell-2}^{(a)}(\nt)=0,\quad P_{\ell-2}^{(a)}(\nt)=0,\quad
 \ell^{(a-1)}<\ell-1,\quad \text{$\ell=\lb^{(a)}=\ell^{(a)}+1$ finite.}
\end{equation*}
We show that these conditions cannot be met simultaneously. Let $p<\ell-1$
be maximal such that $m_p^{(a)}(\nt)>0$; if no such $p$ exists set $p=0$.
By \eqref{eq:convex}, $P_{\ell-2}^{(a)}(\nt)=0$ is only possible
if $P_i^{(a)}(\nt)=0$ for $p\le i\le \ell-1$. By \eqref{eq:Pm C}
this requires $m_i^{(a-1)}(\nt)=0$ for $p<i<\ell-1$ so that
$\ell^{(a-1)}<\ell-1$ implies $\ell^{(a-1)}\le p$. But there is
a singular string of length $p$ in $(\nt,\Jt)^{(a)}$ which contradicts
$\ell^{(a)}=\ell-1>p$.

Finally for $a=n$ the only problem occurs if
\begin{equation*}
m_{\ell-2}^{(n)}(\nt)=0,\quad P_{\ell-2}^{(n)}(\nt)=0,\quad
\ell^{(n-1)}<\ell-1,\quad \text{$\ell=\lb^{(n)}$ finite.}
\end{equation*}
By convexity \eqref{eq:convex}, $P_i^{(n)}(\nt)=0$ for $p\le i\le \ell$
where $p<\ell$ is largest such that $m_p^{(n)}(\nt)>0$. Then by \eqref{eq:Pm C}
we also have $m_i^{(n-1)}(\nt)=0$ for $p<i<\ell$ so that $\ell^{(n-1)}<\ell-1$
implies $\ell^{(n-1)}\le p$. But there is a singular string of length $p$
in $(\nt,\Jt)^{(n)}$ which contradicts $\lb^{(n)}=\ell>p$.
\end{proof}

\begin{proof}[Proof of (III) for $C_n^{(1)}$]
Here $b^\natural=\overline{1}$, $\Hb(b\otimes b')=0$ if $b\le
b'$ and $H(b\otimes b')=1$ otherwise.

If $L=1$ then the path is $1$, the rigged configuration is empty,
and both sides of \eqref{eq:cc=D} are zero.

Here \eqref{eq:Deltacc} and \eqref{eq:Delta2cc} are given by
\begin{align}\label{eq:change cc C1}
\Delta(cc(\nt,\Jt))&=\alpha_1^{(1)} \\
\label{eq:dc D C1} \Hb(b_L\otimes b_{L-1})
 &=\chi(\ell^{(1)}=1)
\end{align}
where $\ell^{(i)}$ is determined by the algorithm $\delta$.
Note that there is no contribution from $\lb^{(1)}$ in \eqref{eq:dc D C1}
since $\lb^{(1)}>1$.

Let $\lt^{(a)}$ and $\tilde{\overline{\ell}}^{(a)}$ be the length
of the selected strings defined by the algorithm $\delta$ on
$(\ntt,\Jtt)= \delta(\nt,\Jt)$. Note that if $\ell^{(1)}=1$ then
\eqref{eq:cv C} implies that $b_{L-1}<b_L$ so that both sides of
\eqref{eq:dc D C1} are 1. If $\ell^{(1)}>1$ then
$\tilde{\ell}^{(a)}<\ell^{(a)}$ for $1\le a<n$ so that $b_{L-1}\ge b_L$
and both sides of \eqref{eq:dc D C1} are 0. For $b_L=\overline{n}$,
$\tilde{\overline{\ell}}^{(n)}<\lb^{(n)}$ unless $\lb^{(n)}=2$.
But note that in this case $\ell^{(n)}=1$ and hence $\ell^{(1)}=1$
which contradicts our assumptions. If $\tilde{\overline{\ell}}^{(n)}<
\lb^{(n)}$ then also $\tilde{\overline{\ell}}^{(a)}<\lb^{(a)}$ which
implies that $b_{L-1}\ge b_L$. Hence both sides of \eqref{eq:dc D C1} are 0.

To prove \eqref{eq:change cc C1}, by \eqref{eq:cc} and
\eqref{eq:ch m C1} we have
\begin{equation*}
\begin{split}
&cc(\ntt)=\frac{1}{2}
 \sum_{i,j\ge 1} \sum_{a,b=1}^n \min(t_b i,t_a j)(\alpha_a|\alpha_b)\\
 &\times \left(m_i^{(a)}-\chi(a<n)(\delta_{i,\ell^{(a)}}-\delta_{i,\ell^{(a)}-1})
  -\delta_{i,\frac{\lb^{(a)}}{\f_a}}
  +\delta_{i,\frac{\lb^{(a)}}{\f_a}-1}\right)\\
 &\times \left(m_j^{(b)}-\chi(b<n)(\delta_{j,\ell^{(b)}}-\delta_{j,\ell^{(b)}-1})
  -\delta_{j,\frac{\lb^{(b)}}{\f_b}}
  +\delta_{j,\frac{\lb^{(b)}}{\f_b}-1}\right).
\end{split}
\end{equation*}
Applying the data for $C_n$ and using \eqref{eq:cv C}, a tedious
but straightforward calculation yields
\begin{multline*}
\Delta cc(\nt)=\sum_{a=1}^n \sum_{i\ge 1}
 \left(P_{\f_a i}^{(a)}(\ntt)-P_{\f_a i}^{(a)}(\nt)\right)
 \left(m_i^{(a)}-\chi(a<n)\delta_{i,\ell^{(a)}}-\delta_{\f_a i,\lb^{(a)}}\right)\\
 +\sum_{i\ge 1}m_i^{(1)}.
\end{multline*}
For $\Delta |\Jt|$ we obtain from the algorithm $\delta'$
\begin{equation*}
\Delta |\Jt|=\sum_{a=1}^n \sum_{i\ge 1}
 \left(P_{\f_a i}^{(a)}(\nt)-P_{\f_a i}^{(a)}(\ntt)\right)
 \left(m_i^{(a)}-\chi(a<n)\delta_{i,\ell^{(a)}}-\delta_{\f_a i,\lb^{(a)}}\right).
\end{equation*}
Hence altogether, using $\sum_{i\ge 1}m_i^{(1)}=\alpha_1^{(1)}$,
we obtain \eqref{eq:change cc C1}.
\end{proof}

\subsection{Proof for type $A_{2n}^{(2)}$}

The proofs of (I) and (II) are analogous to the previous cases. In
particular, the proof of (II) is very similar to that for type
$C_n^{(1)}$.

\begin{proof}[Proof of (III) for $A_{2n}^{(2)}$]
Here $b^\natural=\phi$, $\Hb(b\otimes b')=0$ if $b\le b'$,
$\Hb(b\otimes b')=2$ if $b>b'$ or $b\otimes b'=\phi\otimes \phi$
and $H(b\otimes b')=0$ otherwise.

If $L=1$ then the path is $1$ or $\phi$. In the former case,
the rigged configuration is empty, and both sides of \eqref{eq:cc=D} are zero.
In the other case both sides of \eqref{eq:cc=D} are 1.

Here \eqref{eq:Deltacc} and \eqref{eq:Delta2cc} are given by
\begin{align}\label{eq:change cc A2e}
\Delta(cc(\nt,\Jt))&=2\alpha_1^{(1)}-\chi(\ell^{(n)}=1)\\
\label{eq:dc D A2e} \Hb(b_L\otimes b_{L-1})
 &=2\chi(\ell^{(1)}=1)-\chi(\ell^{(n)}=1)+\chi(\lt^{(n)}=1)
\end{align}
where $\ell^{(i)}$ is determined by the algorithm $\delta$.
Note that there is no contribution from $\lb^{(1)}$ in \eqref{eq:dc D A2e}
since $\lb^{(1)}>1$.

Equation \eqref{eq:dc D A2e} can be checked in a similar fashion as
to the other cases.

To prove \eqref{eq:change cc A2e},
applying the data for $B_n$ and using \eqref{eq:cv C}, a tedious
but straightforward calculation yields
\begin{multline*}
\Delta cc(\nt)=2\sum_{a=1}^n \sum_{i\ge 1}
 \left(P_i^{(a)}(\ntt)-P_i^{(a)}(\nt)\right)
 \left(m_i^{(a)}-\chi(a<n)\delta_{i,\ell^{(a)}}-\delta_{i,\lb^{(a)}}\right)\\
 -\chi(\ell^{(n)}=1)+2\sum_{i\ge 1}m_i^{(1)}.
\end{multline*}
For $\Delta |\Jt|$ we obtain from the algorithm $\delta'$
\begin{equation*}
\Delta |\Jt|=2\sum_{a=1}^n \sum_{i\ge 1}
 \left(P_i^{(a)}(\nt)-P_i^{(a)}(\ntt)\right)
 \left(m_i^{(a)}-\chi(a<n)\delta_{i,\ell^{(a)}}-\delta_{i,\lb^{(a)}}\right).
\end{equation*}
Hence altogether, using $\sum_{i\ge 1}m_i^{(1)}=\alpha_1^{(1)}$,
we obtain \eqref{eq:change cc A2e}.
\end{proof}

\subsection{Proof for type $A_{2n-1}^{(2)}$}

\begin{proof}[Proof of (I) for $A_{2n-1}^{(2)}$]
The proof that $\rho$ is dominant for $b=k$ with $1\le k\le n$ is
analogous to the other types. For $b=\overline{k}$ with $1\le k\le
n$, $\rho$ is not dominant if $\la_k=\la_{k-1}$. In this case
$P_i^{(k-1)}(\nt)=0$ for $i\ge \ell$ where $\ell$ is the largest
part of $\nu^{(k-1)}$ by Lemma \ref{lem:asym} and
\eqref{eq:convex}. By \eqref{eq:Pm A2o}, $m_i^{(k)}(\nt)=0$ for
$i>\ell$ so that $\lb^{(k)}\le \ell$. But since
$P_\ell^{(k-1)}(\nt)=0$ and $m_\ell^{(k-1)}>0$, there is a
singular string of length $\ell$ in $(\nt,\Jt)^{(k-1)}$. Hence
$\lb^{(k-1)}\le \ell$ (which contradicts $\lb^{(k-1)}=\infty$
since $\overline{k}=\rk(\nt,\Jt)$) unless
$\ell^{(k-1)}=\lb^{(k)}=\ell$ and $m_\ell^{(k-1)}(\nt)=1$. Since
$m_\ell^{(k)}(\nt)\ge 2$ for $1\le k<n$ and $m_\ell^{(n)}(\nt)\ge
1$, \eqref{eq:Pm A2o} for $a=k-1$ and $i=\ell$ implies that
$m_\ell^{(k-2)}(\nt)=0$ and $P_{\ell-1}^{(k-1)}(\nt)=0$. Hence
$\ell^{(k-2)}<\ell$ and $m_{\ell-1}^{(k-1)}(\nt)=0$ since
otherwise $\ell^{(k-1)}\le \ell-1$. By induction on
$i=\ell-1,\ell-2,\ldots,1$ \eqref{eq:Pm A2o} for $a=k-1$ implies
that $m_i^{(k-2)}(\nt)=0$ and $P_{i-1}^{(k-1)}(\nt)=0$ which in
turn requires $m_{i-1}^{(k-1)}(\nt)=0$ since else $\ell^{(k-1)}\le
i-1$. But then $m_i^{(k-2)}(\nt)=0$ for all $1\le i\le \ell$ so
that $\ell^{(k-2)}>\ell$ which contradicts our assumptions.
\end{proof}

\begin{proof}[Proof of (II) for $A_{2n-1}^{(2)}$]
To prove that $(\ntt,\Jtt)$ is admissible, one finds similarly to the
proof of type $D_n^{(1)}$ that the only problem occurs if
\begin{multline*}
m_{\ell-1}^{(a)}=0,\quad P_{\ell-1}^{(a)}(\nt)=0,\quad
\text{$\ell^{(a-1)}<\ell$ (resp. $\lb^{(a+1)}<\ell$ for $1\le a<n$)},\\
\text{$\ell$ finite},\quad \text{$\ell=\ell^{(a)}$ (resp. $\ell=\lb^{(a)}$
for $1\le a<n$)}.
\end{multline*}
Analogous to the case $D_n^{(1)}$ it can be shown that these conditions
cannot hold simultaneously.
\end{proof}

\begin{proof}[Proof of (III) for $A_{2n-1}^{(2)}$]
Here $b^\natural=\overline{1}$, $\Hb(\overline{1}\otimes 1)=2$,
$\Hb(b\otimes b')=0$ if $b\le b'$ and $H(b\otimes b')=1$ otherwise.

If $L=1$ then the path is $1$, the rigged configuration is empty,
and both sides of \eqref{eq:cc=D} are zero.

Here \eqref{eq:Deltacc} and \eqref{eq:Delta2cc} are given by
\begin{align}\label{eq:change cc A2o}
\Delta(cc(\nt,\Jt))&=\alpha_1^{(1)} \\
\label{eq:dc D A2o} \Hb(b_L\otimes b_{L-1})
 &=\chi(\ell^{(1)}=1)+\chi(\lb^{(1)}=1)
\end{align}
where $\ell^{(i)}$ is determined by the algorithm $\delta$.

The proof that \eqref{eq:dc D A2o} holds is very similar to the previous
cases.

To prove \eqref{eq:change cc A2o} we apply the data for $C_n$ to
\eqref{eq:cc}. Using \eqref{eq:ch m A2o} and \eqref{eq:cv A2o}
a tedious but straightforward calculation yields
\begin{multline*}
\Delta cc(\nt)=\sum_{a=1}^n \sum_{i\ge 1}
 t_a^\vee \left(P_i^{(a)}(\ntt)-P_i^{(a)}(\nt)\right)
 \left(m_i^{(a)}-\delta_{i,\ell^{(a)}}-\chi(a<n)\delta_{i,\lb^{(a)}}\right)\\
 +\sum_{i\ge 1}m_i^{(1)}.
\end{multline*}
For $\Delta |\Jt|$ we obtain from the algorithm $\delta'$
\begin{equation*}
\Delta |\Jt|=\sum_{a=1}^n \sum_{i\ge 1}
 t_a^\vee \left(P_i^{(a)}(\nt)-P_i^{(a)}(\ntt)\right)
 \left(m_i^{(a)}-\delta_{i,\ell^{(a)}}-\chi(a<n)\delta_{i,\lb^{(a)}}\right).
\end{equation*}
Hence altogether, using $\sum_{i\ge 1}m_i^{(1)}=\alpha_1^{(1)}$,
we obtain \eqref{eq:change cc A2o}.
\end{proof}

\subsection{Proof for type $D_{n+1}^{(2)}$}

\begin{proof}[Proof of (I) for $D_{n+1}^{(2)}$]
The proof proceeds as before
except in the cases $b=n$ and $b=0$ (there is nothing to prove for
$b=\phi$). Suppose $b=n$ and $\rho$ is not dominant. Since $\la$
is dominant, $\la_n=0$. Then it can be deduced that
$P_i^{(n)}(\nt)=0$ for $i\ge \ell$ where $\ell$ is the largest
part in $\nu^{(n)}$ by Lemma \ref{lem:asym} and \eqref{eq:convex}
and the admissibility of $\nt$. By \eqref{eq:Pm D2} it follows
that $m_i^{(n-1)}(\nt)=0$ for $i>\ell$ so that $\ell^{(n-1)}\le
\ell$. But there is a singular string of length $\ell$ in
$(\nt,\Jt)^{(n)}$ since $P_\ell^{(n)}(\nt)=0$ which contradicts
$\ell^{(n)}=\infty$. To prove that $b=0$ cannot occur if $\la_n=0$
we find as for the case $k=n$ that $\ell^{(n-1)}\le \ell$ and that
there is a singular string of length $\ell$ in $(\nt,\Jt)^{(n)}$
since $P_\ell^{(n)}(\nt)=0$. For the $b=0$ case (Q) must hold so
that there must be a quasisingular string of length
$\ell^{(n)}<\ell$ in $(\nt,\Jt)^{(n)}$. But observe that there is
a singular string of length $\ell>\ell^{(n)}$ in $(\nt,\Jt)^{(n)}$
which contradicts $\lb^{(n)}=\infty$.
\end{proof}

\begin{proof}[Proof of (II) for $D_{n+1}^{(2)}$]
Next we need to show that $(\ntt,\Jtt)\in\RC(\rho,\mut)$. The case
that $(\ntt,\Jtt)^{(a)}$ is admissible for $1\le a<n$ works as
usual. Consider $a=n$. First note that there is no problem in case
(Q,S) setting the new string of length $\lb^{(n)}-1$ to be
quasisingular since the string of length $\lb^{(n)}-1$ is not
singular by definition so that $P_{\lb^{(n)}-1}^{(n)}(\nt)>0$ and
also $P_{\lb^{(n)}-1}^{(n)}(\ntt)>0$ by \eqref{eq:cv D2}. The only
problem occurs if
\begin{equation*}
m_{\ell-1}^{(n)}(\nt)=0,\quad P_{\ell-1}^{(n)}(\nt)=\text{0 or 1},\quad
\ell^{(n-1)}<\ell,\quad \text{$\ell=\ell^{(n)}$ finite.}
\end{equation*}
Note that $P_i^{(n)}(\nt)$ is always even so that $P_{\ell-1}^{(n)}=1$
is impossible. The proof that these conditions cannot hold simultaneously
works as usual.
\end{proof}

\begin{proof}[Proof of (III) for $D_{n+1}^{(2)}$]
Here $b^\natural=\phi$ and $\Hb(\phi\otimes\phi)=2$,
$\Hb(b\otimes \phi) =\Hb(\phi\otimes b)=1$ if $b\neq \phi$,
$\Hb(b\otimes b')=0$ if $b,b'\neq \phi$, $b\le b'$ and $b\neq b'$ if $b=0$,
and $\Hb(b\otimes b')=2$ if $b>b'$ or $b=b'=0$.

If $L=1$ then the path is either $1$ or $\phi$. In the former
case the rigged configuration is empty, and both sides of
\eqref{eq:cc=D} are zero. In the latter case it is also not hard
to check that both sides of \eqref{eq:cc=D} are 1.

Here \eqref{eq:Deltacc} and \eqref{eq:Delta2cc} are given by
\begin{align}\label{eq:change cc D2}
\Delta(cc(\nt,\Jt))&=2 \alpha_1^{(1)} -\chi(\ell^{(n)}=1) \\
\label{eq:dc D D2} \Hb(b_L\otimes b_{L-1})
 &=2\chi(\ell^{(1)}=1)-\chi(\ell^{(n)}=1) +\chi(\lt^{(n)}=1)
\end{align}
where $\ell^{(i)}$ and $\lb^{(i)}$ are determined by the algorithm
$\delta$. To obtain \eqref{eq:dc D D2} we used the fact that by
definition $\lb^{(1)}>1$. Here $\lt^{(a)}$ is defined by the
algorithm $\delta$ on $(\ntt,\Jtt)=\delta'(\nt,\Jt)$.

It can be checked directly that \eqref{eq:dc D D2} holds. For
example, if $\ell^{(1)}=1$ it follows that $\lt^{(a)}\ge
\ell^{(a+1)}$ for $1\le a<n$. Hence if $b_L\le n$ then
$b_{L-1}<b_L$ and both sides of \eqref{eq:dc D D2} yield 2. If
$b_L=\phi$ then both sides of \eqref{eq:dc D D2} are 1 for
$b_{L-1}\neq \phi$ and 2 if $b_{L-1}=\phi$. If $b_L=0$
then both sides of \eqref{eq:dc D D2} are 2 if $b_{L-1}\le 0$.
Note that $b_{L-1}=\phi$ or $b_{L-1}\ge \overline{n}$ is not
possible. Finally, if $b_L\ge \overline{n}$ then
$\tilde{\overline{\ell}}^{(a)}\ge \lb^{(a)}$ and $b_{L-1}<b_L$.
Note that $b_{L-1}=\phi$ is not possible in this case since
$\ell^{(n)}>1$ which implies $\lt^{(n)}>1$. Both sides of
\eqref{eq:dc D D2} yield 2 in this case. If $\ell^{(1)}>1$ then
$\lt^{(a)}<\ell^{(a)}$ for $1\le a\le n$ and the cases can be
checked in a similar fashion as before.

To prove \eqref{eq:change cc D2}, from \eqref{eq:cc} and
\eqref{eq:ch m C1} we obtain
\begin{equation*}
\begin{split}
cc(\ntt)=&\frac{1}{2}
 \sum_{i,j\ge 1} \sum_{a,b=1}^n \min(i,j) (\alpha_a|\alpha_b)\\
 &\times (m_i^{(a)}-\delta_{i,\ell^{(a)}}-\delta_{i,\lb^{(a)}}
 +\delta_{i,\ell^{(a)}-1}+\delta_{i,\lb^{(a)}-1})\\
 &\times(m_j^{(b)}-\delta_{j,\ell^{(b)}}-\delta_{j,\lb^{(b)}}
 +\delta_{j,\ell^{(b)}-1}+\delta_{j,\lb^{(b)}-1}).
\end{split}
\end{equation*}
Expanding out and using \eqref{eq:cv D2} a tedious but straightforward
calculation yields
\begin{multline*}
\Delta cc(\nt)=\sum_{a=1}^n \sum_{i\ge 1}
 t_a^\vee \left(P_i^{(a)}(\ntt)-P_i^{(a)}(\nt)\right)
 \left( m_i^{(a)}-\delta_{i,\ell^{(a)}}-\delta_{i,\lb^{(a)}}\right)\\
+2\sum_{i\ge 1}m_i^{(1)}-\chi(\lb^{(n)}=\infty)\chi(\ell^{(n)}<\infty).\\
\end{multline*}
For $\Delta |\Jt|$ we obtain from the algorithm $\delta'$
\begin{multline*}
\Delta |\Jt|=\sum_{a=1}^n \sum_{i\ge 1}
 t_a^\vee\left(P_i^{(a)}(\nt)-P_i^{(a)}(\ntt)\right)
 \left(m_i^{(a)}-\delta_{i,\ell^{(a)}}-\delta_{i,\lb^{(a)}}\right)\\
+\chi(\lb^{(n)}=\infty)\chi(1<\ell^{(n)}<\infty)
\end{multline*}
where the last term comes from the fact that in case (Q) a
quasisingular string is changed into a singular string. Hence
altogether, using $\sum_{i\ge 1}m_i^{(1)}=\alpha_1^{(1)}$ and the
fact that $\lb^{(n)}=\infty$ if $\ell^{(n)}=1$ by the algorithm
$\delta$, we obtain \eqref{eq:change cc D2}.
\end{proof}

\subsection{Proof for type $A_{2n}^{(2)\dagger}$}

\begin{proof}[Proof of (I) for $A_{2n}^{(2)\dagger}$]
The only case that proceeds differently than before is
$b=0$. Suppose $\la_n=0$. Let $\ell$ be the longest part of
$\nu^{(n)}$. As in the proof of the $D^{(2)}_{n+1}$ case,
$P^{(n)}_\ell(\nt)=0$ where $\ell \ge \ell^{(n-1)}$. If $\ell$ is
odd then this is a contradiction of the admissibility of $\nt$;
see \eqref{eq:a2doddvac}. If $\ell$ is even then $J^{(n,\ell)}$ is
singular and $\ell^{(n)}<\ell$ (as $\ell^{(n)}$ is odd and $\ell$
is the longest part), contradicting $b=0$.
\end{proof}

\begin{proof}[Proof of (II) for $A_{2n}^{(2)\dagger}$]
The admissibility of $(\ntt,\Jtt)$ for $1\le a<n$ is as before.
Let $a=n$. We first observe that in all cases,
\begin{equation}
\ell^{(1)}\le\ell^{(2)}\le\dotsm\le\ell^{(n)}\le
\lb^{(n)}\le\lb^{(n-1)}\le\dotsm\le\lb^{(1)},
\end{equation}
with $\ell^{(n)}$ odd and $\lb^{(n)}$ even (when they are finite).
We also note that by \eqref{eq:cv C},
\begin{equation} \label{eq:onedrop}
P^{(n)}_i(\ntt) \ge P^{(n)}_i(\nt) - 1
\end{equation}
with equality if and only if $\ell^{(n-1)} \le i<\ell^{(n)}$.

Let us verify \eqref{eq:a2doddvac} for $(\ntt,\Jtt)$. Let $i$ be
odd such that $m_i^{(n)}(\ntt)>0$. Suppose first that
$m_i^{(n)}(\nt)>0$. By \eqref{eq:a2doddvac} for $(\nt,\Jt)$,
$P^{(n)}_i(\nt)>0$. By \eqref{eq:onedrop} we may assume that
$P^{(n)}_i(\nt)=1$ and $\ell^{(n-1)}\le i<\ell^{(n)}$. But then
$J^{(n,i)}$ was quasisingular, which is a contradiction to the
definition of $\delta$. So suppose $m_i^{(n)}(\nt)=0$. Since
$m_i^{(n)}(\ntt)>0$ we are in case $(Q,S)$ with $i=\lb^{(n)}-1$.
In case $(Q,S)$ $\ell^{(n)}<\lb^{(n)}$, so $\ell^{(n)} \le i <
\lb^{(n)}$. Now $i\not=\ell^{(n)}$ since $m_i^{(n)}(\nt)=0$. So
$\ell^{(n)} < i < \lb^{(n)}$ with $\ell^{(n)}$ and $i$ odd. By
\eqref{eq:cv C} $P^{(n)}_i(\ntt)=P^{(n)}_i(\nt)$. There is only
a problem if $P^{(n)}_i(\nt)=0$. By \eqref{eq:convex} it follows
that $P^{(n)}_{i-1}(\nt)=P^{(n)}_{i+1}(\nt)=0$. Since $i-1$ is
even, if $m_{i-1}^{(n)}(\nt)>0$ then $J^{(n,i-1)}$ would have been
singular with $\ell^{(n)}<i-1<\lb^{(n)}$, contradicting the choice
of $\lb^{(n)}$. So $m_{i-1}^{(n)}(\nt)=0$. Applying
\eqref{eq:convex} again, $P_{i-2}^{(n)}(\nt)=0$. Since
$(\nt,\Jt)$ was admissible and $i-2$ is odd, by
\eqref{eq:a2doddvac} it follows that $m_{i-2}^{(n)}(\nt)=0$.
Continuing in this manner, a contradiction is reached since
$P^{(n)}_{\ell^{(n)}}(\nt) > 0$.

Now suppose $i$ is even. It must be checked that
$P^{(n)}_i(\ntt)\ge0$. The only problem is if
$P^{(n)}_i(\nt)=0$ and $\ell^{(n-1)}\le i<\ell^{(n)}$. If
$m_i^{(n)}(\nt)>0$ then $\delta$ would have chosen the singular
partition $J^{(n,i)}$. So $m_i^{(n)}(\nt)=0$. By \eqref{eq:convex}
it follows that $P^{(n)}_{i+1}(\nu)=0$. Arguing as above but with
the index increasing from $i$, a contradiction is reached since
$P^{(n)}_{\ell^{(n)}}(\nt)>0$.
\end{proof}

\begin{proof}[Proof of (III) for $A_{2n}^{(2)\dagger}$]
One has $b^\natural=1$, $\Hb(b_2'\otimes b_1')=0$ if
$b_2'\le b_1'$ (except for $\Hb(0\otimes 0)=1$), and
$\Hb(b'_2\otimes b'_1)=1$ for $b'_2>b'_1$.

If $L=1$ then the path is $1$, the rigged configuration is empty,
and both sides of \eqref{eq:cc=D} are zero.

Here \eqref{eq:Deltacc} and \eqref{eq:Delta2cc} are given by
\begin{align}\label{eq:change cc A2d}
\Delta(cc(\nt,\Jt))&=\alpha_1^{(1)} \\
\label{eq:dc D A2d} \Hb(b_L\otimes b_{L-1})
 &=\chi(\ell^{(1)}=1)
\end{align}
where $\ell^{(i)}$ and $\lb^{(i)}$ are determined by the algorithm
$\delta$. The term $\chi(\lb^{(1)}=1)$ disappears since the
definition of the algorithm forces $\lb^{(1)}\ge 2$. The proof of
\eqref{eq:change cc A2d} is very similar to that in the
$D_{n+1}^{(2)}$ case.

Straightforward computations yield
\begin{equation*}
\begin{split}
  \Delta cc(\nt) &=
  \sum_{a,b,i} (\alpha_a|\alpha_b)
  (\chi(i\ge\ell^{(b)})+\chi(i\ge\lb^{(b)}))
  (m_i^{(a)}(\nt)-\delta_{i,\ell^{(a)}}-\delta_{i,\lb^{(a)}}) \\
  &+ \chi(\ell^{(1)}<\infty) + \chi(\lb^{(1)}<\infty) -
  \frac{1}{2} \chi(\ell^{(n)}<\infty)+\frac{1}{2}
  \chi(\lb^{(n)}<\infty)
\end{split}
\end{equation*}
and
\begin{equation*}
\begin{split}
& \sum_{a,i} (P_i^{(a)}(\ntt)-P_i^{(a)}(\nt))
(m_i^{(a)}(\nt)-\delta_{i,\ell^{(a)}}-\delta_{i,\lb^{(a)}}) \\
&=\sum_{a,b,i} (\alpha_a|\alpha_b)
  (\chi(i\ge\ell^{(b)})+\chi(i\ge\lb^{(b)}))
  (m_i^{(a)}(\nt)-\delta_{i,\ell^{(a)}}-\delta_{i,\lb^{(a)}}) \\
  &- \sum_i m_i^{(a)}(\nt) + \chi(\ell^{(1)}<\infty) +
    \chi(\lb^{(1)}<\infty).
\end{split}
\end{equation*}
Together these yield
\begin{equation*}
\begin{split}
\Delta cc(\nt)&=\sum_{a,i} (P_i^{(a)}(\ntt)-P_i^{(a)}(\nt))
(m_i^{(a)}(\nt)-\delta_{i,\ell^{(a)}}-\delta_{i,\lb^{(a)}}) \\
&+ \alpha^{(1)}_1 - \frac{1}{2} \chi(\ell^{(n)}<\infty)
+\frac{1}{2} \chi(\lb^{(n)}<\infty).
\end{split}
\end{equation*}
One can also show that
\begin{equation*}
\begin{split}
  \Delta |\Jt| &= \sum_{a,i}
  (P_i^{(a)}(\nt)-P_i^{(a)}(\ntt))(m_i^{(a)}(\nt)-\delta_{i,\ell^{(a)}}-\delta_{i,\lb^{(a)}})
  \\
  &+ \frac{1}{2} \chi(\ell^{(n)}<\infty) -
  \frac{1}{2}
  \chi(\lb^{(n)}<\infty).
\end{split}
\end{equation*}
This proves \eqref{eq:change cc A2d}.
\end{proof}

\end{document}